\documentclass[12pt]{article}
\usepackage{multicol}
\usepackage{amsmath, amssymb, amscd, enumerate, amsfonts, latexsym, amsxtra, 
amsthm, bm, multirow}
\usepackage{indentfirst}
\newtheorem{The}{Theorem}[section]
\newtheorem{Prop}[The]{Proposition}
\newtheorem{Lem}[The]{Lemma}
\newtheorem{Cor}[The]{Corollary}

\begin{document}
\allowdisplaybreaks[3]
\centerline{\Large Certain linear isomorphisms for  hyperalgebras} \vspace{3mm}
\centerline{\Large relative to  a Chevalley group} \vspace{7mm}
\centerline{Yutaka Yoshii 
\footnote{ E-mail address: yutaka.yoshii.6174@vc.ibaraki.ac.jp}}  \vspace{5mm}
\centerline{College of Education,   
Ibaraki University,}
\centerline{2-1-1 Bunkyo, Mito, Ibaraki, 310-8512, Japan}
\begin{abstract}
Let $G$ be a simply connected and simple algebraic group defined and split over 
a finite prime field $\mathbb{F}_p$ of $p$ elements. In this paper, using an 
$\mathbb{F}_p$-linear map splitting Frobenius endomorphism on a hyperalgebra 
relative to $G$, we obtain some $\mathbb{F}_p$-linear isomorphisms induced by 
multiplication in the hyperalgebra. 
\end{abstract}

\section{Introduction} 
Let $\mathbb{F}_p$ be a finite prime field of $p$ elements and 
$\overline{\mathbb{F}}_p$ its algebraic closure. Let $G$ be a simply connected and 
simple algebraic $\overline{\mathbb{F}}_p$-group defined and split over $\mathbb{F}_p$. Let 
$B^+$ and $B$ be the Borel subgroups of $G$ corresponding to the positive 
and the negative roots respectively. Let $U^+$ and $U$ be the unipotent radicals of 
$B^+$ and $B$ respectively.  

There is a unique $\overline{\mathbb{F}}_p$-algebra endomorphism 
${\rm Fr}: {\rm Dist}(G) \rightarrow {\rm Dist}(G)$ called the Frobenius endomorphism 
on the hyperalgebra ${\rm Dist}(G)$, which corresponds to the geometric  
Frobenius map $G \rightarrow G$. Gros and Kaneda introduced a linear map  
${\rm Dist}(G) \rightarrow {\rm Dist}(G)$ splitting the   
Frobenius endomorphism in \cite{gros-kaneda11}, which will be denoted by 
${\rm Fr}'$ in this paper. This map, which is originally based on 
a map defined on a quantum group by Lusztig \cite[Lemma 8.6]{lusztig90}, has a lot of  
interesting properties.  

The author is interested in the map and expects its various applications. For example, 
in \cite{yoshii17}, he constructed the primitive idempotents in the hyperalgebra 
${\rm Dist}(G_r)$ of the 
$r$-th Frobenius kernel $G_r$ for $G={\rm SL}_2(\overline{\mathbb{F}}_p)$ whose sum is 
the unity $1$, using the Frobenius splitting ${\rm Fr}'$. Moreover, it turns out that these 
primitive idempotents have further applications (see \cite{yoshii18} and \cite{yoshii22-1}). 

By the way, in arguing construction of  
the above primitive idempotents of 
${\rm Dist}(G_r)$ for $G={\rm SL}_2(\overline{\mathbb{F}}_p)$, the author uses 
the fact that the map 
\[ {\rm Dist}(G_n) \otimes_{\overline{\mathbb{F}}_p} {\rm Fr}'^n 
\left( {\rm Dist}(G_{n'-n})\right) \rightarrow 
{\rm Dist}(G_{n'})\]
for $n,n' \in \mathbb{Z}_{>0}$ with $n'>n$ 
induced by multiplication is an isomorphism as $\overline{\mathbb{F}}_p$-vector spaces 
(see \cite[Proposition 2.3]{yoshii17}). 
So if  bijectivity of the map holds for general $G$, further application of the map 
to representation theory for $G$ can be expected. 

In this paper, for a general simply connected and simple algebraic group $G$ and its primary 
subgroups, we claim that several maps of this kind induced by multiplication 
give  various linear isomorphisms. 
Of course, the argument is much more complicated than the simplest case 
$G={\rm SL}_2(\overline{\mathbb{F}}_p)$. For convenience, we would rather use 
the base change $G_{\mathbb{F}_p}$ and its hyperalgebra 
${\rm Dist}(G_{\mathbb{F}_p})$  than $G$ and ${\rm Dist}(G)$ 
themselves. The results over $\mathbb{F}_p$ immediately imply the corresponding ones over 
$\overline{\mathbb{F}}_p$. 

Main results will be given in Sections 4 and 5.  In Section 3, 
we give several commutation formulas as a preparation 
for proving main results. 
In Section 4, we give several 
$\mathbb{F}_p$-linear isomorphisms induced by multiplication 
on ${\rm Dist}(U_{\mathbb{F}_p}^+)$. First we give a result on hyperalgebras with respect to 
 $U_{r, \mathbb{F}_p}^+$ and then give one with respect to 
$U_{\mathbb{F}_p}^+$. We also use the fact that the restriction of ${\rm Fr}'$ to 
${\rm Dist}(U_{\mathbb{F}_p}^+)$ is an $\mathbb{F}_p$-algebra homomorphism which is proved 
by Gros and Kaneda in \cite{gros-kaneda11}. 
Of course, by symmetry, these results for $U^+$ imply the similar ones 
for $U$. Finally, in Section 5, we give several 
$\mathbb{F}_p$-linear isomorphisms induced by multiplication 
on ${\rm Dist}(G_{\mathbb{F}_p})$. These results imply the similar ones for $B^+$ and $B$.

\section{Preliminaries}
Let $p$ be a prime number and  $G$  a simply connected and simple algebraic group 
defined and split over 
$\mathbb{F}_p$.  If $H$ is a connected subgroup of $G$ defined over $\mathbb{F}_p$, 
let $H_r$ be its $r$-th Frobenius kernel for $r \in \mathbb{Z}_{>0}$.

Let $T$ be a split maximal torus of $G$ defined over $\mathbb{F}_p$. 
Let $X(T)$ be the  group of rational characters of $T$. Let $\Phi$ be the root system 
relative to the pair $(G,T)$. Let $\Delta$ be the set of simple roots and $\Phi^+$ (resp.  $\Phi^-$)  the set of positive 
(resp. negative) roots. 
Let $B$ (resp. $B^+$) be the Borel subgroup of $G$ corresponding to 
$\Phi^-$ (resp. $\Phi^+$). Let $U$ (resp. $U^+$) be the unipotent radical of 
$B$ (resp. $B^+$). Let $W=N_{G}(T)/T$ be the Weyl group. Then the action of $N_{G}(T)$ 
on $T$ 
by conjugation induces that of $W$  on $X(T)$, and hence on the 
euclidean space $\mathbb{E}=X(T) \otimes_{\mathbb{Z}} \mathbb{R}$.  
Let $\langle \cdot , \cdot \rangle$ be a $W$-invariant inner product in  $\mathbb{E}$. 
For $\lambda \in \mathbb{E}$, 
we set $||\lambda||=\sqrt{\langle \lambda, \lambda \rangle}$ and call it the length of 
$\lambda$. For $\lambda, \mu \in \mathbb{E}$, the symbol $\angle(\lambda, \mu)$ denotes the angle of $\lambda$ and $\mu$. 
For a root $\alpha \in \Phi$, we set 
$\alpha^{\vee}= 2\alpha/\langle \alpha, \alpha \rangle$ and call it the coroot 
of $\alpha$. In this paper, we assume that $T$ has rank $l\ (\in \mathbb{Z}_{> 0})$ 
and hence  we can write $\Delta = \{ \alpha_1, \dots, \alpha_l\}$. Each $\beta \in \Phi^+$ 
(resp. $\beta \in \Phi^-$) can be uniquely written as  $\beta= \sum_{i=1}^l c_i \alpha_i$ 
with integral coefficients $c_i$ all nonnegative (resp. all nonpositive). 
Then we define ${\rm ht}(\beta)$ as the integer 
$\sum_{i=1}^l c_i$ and call it the height of $\beta$. 
In $X(T)$, we define a partial order $\geq$ as follows:
\[ 
\mbox{$\mu \geq \lambda$ if $\mu - \lambda \in \sum_{i=1}^{l} 
\mathbb{Z}_{\geq 0} \alpha_i$.} 
\]
For a root 
$\alpha \in \Phi$, let $s_{\alpha}$ be the reflection for the hyperplane 
orthogonal to $\alpha$. Then $W$ is generated by all $s_{\alpha}$ for $\alpha \in \Delta$.  

Let $\mathfrak{g}_{\mathbb{C}}$ be the simple complex Lie algebra which has $\Phi$ 
as a root system. Let $[\cdot , \cdot]$ be the Lie bracket in $\mathfrak{g}_{\mathbb{C}}$ 
and $\{ e_{\alpha}, h_i\ |\ \alpha \in \Phi, 
1 \leq i \leq l\}$  a Chevalley basis of $\mathfrak{g}_{\mathbb{C}}$ which 
satisfies $[e_{\alpha_i}, e_{-\alpha_i}]=h_i$ for $1 \leq i \leq l$. For $\alpha \in \Phi$, 
set $h_{\alpha}=[e_{\alpha}, e_{-\alpha}]$. It is an integral linear combination of the 
$h_i$. Of course we have $h_i=h_{\alpha_i}$. 
Let $\mathcal{U}(\mathfrak{g}_{\mathbb{C}})$ be 
the universal enveloping algebra of $\mathfrak{g}_{\mathbb{C}}$. 
Then $[x,y]=xy-yx$ in $\mathcal{U}(\mathfrak{g}_{\mathbb{C}})$ for $x, y \in \mathfrak{g}_{\mathbb{C}}$. 
For $\alpha \in \Phi$ and $n \in \mathbb{Z}_{\geq 0}$, set 
$e_{\alpha}^{(n)}=e_{\alpha}^n/n!$ in $\mathcal{U}(\mathfrak{g}_{\mathbb{C}})$. 
For convenience, set $e_{\alpha}^{(n)}=0$ for $n \in \mathbb{Z}_{<0}$. 
Let $\mathfrak{g}_{\mathbb{Z}}$ be a $\mathbb{Z}$-span of the Chevalley basis of 
$\mathfrak{g}_{\mathbb{C}}$. 
We define Kostant's 
$\mathbb{Z}$-form $\mathcal{U}_{\mathbb{Z}}$ as a subring of 
$\mathcal{U}(\mathfrak{g}_{\mathbb{C}})$ generated by all 
$e_{\alpha}^{(n)}$ for $\alpha \in \Phi$ and $n \in \mathbb{Z}_{\geq 0}$. 
Clearly $\mathcal{U}_{\mathbb{Z}}$ contains $\mathfrak{g}_{\mathbb{Z}}$. 
For $\alpha \in \Phi$, $c \in \mathbb{Z}$, and $n \in \mathbb{Z}_{\geq 0}$, set 
\[{h_{\alpha} +c \choose n}=\dfrac{\prod_{j=1}^n (h_{\alpha}+c-n+j)}{n!}\]
in $\mathcal{U}(\mathfrak{g}_{\mathbb{C}})$, which lies in $\mathcal{U}_{\mathbb{Z}}$.

For a moment, suppose  that $\Phi$ is of type ${\rm G}_2$ with 
$||\alpha_1|| < ||\alpha_2||$. Then we have 
\[ \Phi^+=\{\alpha_1, \alpha_2, \alpha_1+\alpha_2, 2\alpha_1+\alpha_2, 
3\alpha_1+\alpha_2, 3\alpha_1+2\alpha_2\}. \]
In this case, for simplicity we shall denote 
$e_{\alpha_1}$, $e_{\alpha_2}$, $e_{\alpha_1+\alpha_2}$,  $e_{2\alpha_1+\alpha_2}$, 
$e_{3\alpha_1+\alpha_2}$, $e_{3\alpha_1+2\alpha_2}$ by 
$e_1$, $e_2$, $e_{12}$, $e_{112}$, $e_{1112}$, $e_{11122}$ respectively. 
Without loss of generality, we may assume in this paper 
that these elements are chosen such that 
$[e_1,e_2]=e_{12}$, $[e_1,e_{12}]=2e_{112}$, $[e_1,e_{112}]=3e_{1112}$, 
$[e_2,e_{1112}]=e_{11122}$, and $[e_{112},e_{12}]=3e_{11122}$
in $\mathfrak{g}_{\mathbb{Z}}$ (for example, see \cite[\S 5]{lusztig90}).  

The following formulas are well-known. \\

\begin{Prop}\label{basicformulas}
Let $\alpha, \beta \in \Phi$, $c \in \mathbb{Z}$, and 
$m,n \in \mathbb{Z}_{\geq 0}$. 
In $\mathcal{U}_{\mathbb{Z}}$, 
 the following equalities hold. \\ 

\noindent {\rm (i)} 
$\displaystyle{e_{\alpha}^{(m)} e_{\alpha}^{(n)} = {m+n \choose n} e_{\alpha}^{(m+n)}}$. \\

\noindent {\rm (ii)} $\displaystyle{e_{\alpha}^{(m)} e_{-\alpha}^{(n)}= 
\sum_{k=0}^{{\rm min}\{ m,n\}} e_{-\alpha}^{(n-k)} {h_{\alpha}-m-n+2k \choose k} 
e_{\alpha}^{(m-k)}}$. \\

\noindent {\rm (iii)} 
$\displaystyle{e_{\alpha}^{(m)} {h_{\beta} + c \choose n} = 
{h_{\beta} + c -\langle \alpha, \beta^{\vee} \rangle m \choose n}e_{\alpha}^{(m)}}$. \\

\noindent {\rm (iv)} $e_{\alpha}^{(m)} e_{\beta}^{(n)}= 
e_{\beta}^{(n)}e_{\alpha}^{(m)}$ if $\alpha + \beta \not\in \Phi$ and $\beta \neq -\alpha$. \\

\noindent {\rm (v)}  $\displaystyle{{h_{\alpha} \choose m} {h_{\alpha} \choose n}= 
\sum_{k=0}^{{\rm min}\{ m,n \} }{m+n-k \choose n} {n \choose k} {h_{\alpha} \choose m+n-k}}$.
\end{Prop}
\

Let $\mathcal{U}_{\mathbb{Z}}^+$ (resp. $\mathcal{U}_{\mathbb{Z}}^-$, 
$\mathcal{U}_{\mathbb{Z}}^0$) 
be the subring of $\mathcal{U}_{\mathbb{Z}}$ generated by 
$\{ e_{\alpha}^{(n)} \ |\ \alpha \in \Phi^+, n \in \mathbb{Z}_{\geq 0}\}$ 
(resp. $\{ e_{\alpha}^{(n)} \ |\ \alpha \in \Phi^-, n \in \mathbb{Z}_{\geq 0}\}$, 
$\{ {h_i \choose n} \ |\ i \in \{1, \dots, l\}, n \in \mathbb{Z}_{\geq 0}\}$). 
It is well-known that the subring $\mathcal{U}_{\mathbb{Z}}^+$ 
(resp. $\mathcal{U}_{\mathbb{Z}}^-$) is generated by 
$\{ e_{\alpha_i}^{(n)} \ |\ i \in \{1, \dots, l\}, n \in \mathbb{Z}_{\geq 0}\}$ 
(resp. $\{ e_{-\alpha_i}^{(n)} \ |\ i \in \{1, \dots, l\}, n \in \mathbb{Z}_{\geq 0}\}$).  
We have a triangular decomposition  
$\mathcal{U}_{\mathbb{Z}}=
\mathcal{U}_{\mathbb{Z}}^-\mathcal{U}_{\mathbb{Z}}^0\mathcal{U}_{\mathbb{Z}}^+$. 
In other words, the multiplication 
map  
$\mathcal{U}_{\mathbb{Z}}^- \otimes_{\mathbb{Z}} 
\mathcal{U}_{\mathbb{Z}}^0 \otimes_{\mathbb{Z}}
\mathcal{U}_{\mathbb{Z}}^+ \rightarrow \mathcal{U}_{\mathbb{Z}}$ 
is an isomorphism of additive groups. 
We define an $\mathbb{F}_p$-algebra $\mathcal{U}$ as 
$\mathcal{U}= \mathcal{U}_{\mathbb{Z}} \otimes_{\mathbb{Z}} \mathbb{F}_p$. 
We use the same symbols for images in $\mathcal{U}$ of 
the elements of $\mathcal{U}_{\mathbb{Z}}$ (for example, $e_{\alpha}^{(n)}$, 
${h_i+c \choose n}$, and so on). 
Let $\mathcal{U}^+$ (resp. $\mathcal{U}^-$, $\mathcal{U}^0$) 
be the $\mathbb{F}_p$-subalgebra of $\mathcal{U}$ generated by 
$\{ e_{\alpha}^{(n)} \ |\ \alpha \in \Phi^+, n \in \mathbb{Z}_{\geq 0}\}$ 
(resp. $\{ e_{\alpha}^{(n)} \ |\ \alpha \in \Phi^-, n \in \mathbb{Z}_{\geq 0}\}$, 
$\{ {h_i \choose n} \ |\ i \in \{1, \dots, l\}, n \in \mathbb{Z}_{\geq 0}\}$). Of course we have 
$\mathcal{U}^+= \mathcal{U}_{\mathbb{Z}}^+ \otimes_{\mathbb{Z}} \mathbb{F}_p$, 
$\mathcal{U}^-= \mathcal{U}_{\mathbb{Z}}^- \otimes_{\mathbb{Z}} \mathbb{F}_p$, and 
$\mathcal{U}^0= \mathcal{U}_{\mathbb{Z}}^0 \otimes_{\mathbb{Z}} \mathbb{F}_p$.
Let $\mathcal{U}^{\geq 0}$ (resp. $\mathcal{U}^{\leq 0}$) be the 
$\mathbb{F}_p$-subalgebra of $\mathcal{U}$ generated by $\mathcal{U}^+$ and 
$\mathcal{U}^0$ (resp. $\mathcal{U}^-$ and $\mathcal{U}^0$). 
For $r \in \mathbb{Z}_{>0}$, we denote the subset $\{0, 1, \dots, p^r-1\}$ 
of $\mathbb{Z}$ by $\mathcal{N}_r$. 
For $r \in \mathbb{Z}_{>0}$, let $\mathcal{U}_{r}$ 
be the $\mathbb{F}_p$-subalgebra of $\mathcal{U}$ generated by 
$\{ e_{\alpha}^{(n)} \ |\ \alpha \in \Phi, n \in \mathcal{N}_r\}$. 
Moreover, set $\mathcal{U}_r^+=\mathcal{U}^+ \cap \mathcal{U}_{r}$, 
$\mathcal{U}_r^-=\mathcal{U}^- \cap \mathcal{U}_{r}$, 
$\mathcal{U}_r^0=\mathcal{U}^0 \cap \mathcal{U}_{r}$, 
$\mathcal{U}_r^{\geq 0}=\mathcal{U}^{\geq 0} \cap \mathcal{U}_{r}$, and 
$\mathcal{U}_r^{\leq 0}=\mathcal{U}^{\leq 0} \cap \mathcal{U}_{r}$. 
We can identify 
these algebras with hyperalgebras of various group schemes: 
\[\mathcal{U}={\rm Dist}(G_{\mathbb{F}_p}),\ \ \mathcal{U^+}={\rm Dist}(U_{\mathbb{F}_p}^+),
\ \ \mathcal{U}^-={\rm Dist}(U_{\mathbb{F}_p}),\ \ 
\mathcal{U}^0={\rm Dist}(T_{\mathbb{F}_p}), \] 
\[
\mathcal{U}^{\geq 0}={\rm Dist}(B_{\mathbb{F}_p}^+),\ \ 
\mathcal{U}^{\leq 0}={\rm Dist}(B_{\mathbb{F}_p}),\]
\[\mathcal{U}_r={\rm Dist}(G_{r,\mathbb{F}_p}),\ \ 
\mathcal{U}_r^+={\rm Dist}(U_{r,\mathbb{F}_p}^+),
\ \ \mathcal{U}_r^-={\rm Dist}(U_{r,\mathbb{F}_p}),\ \ 
\mathcal{U}_r^0={\rm Dist}(T_{r,\mathbb{F}_p}), \] 
\[
\mathcal{U}_r^{\geq 0}={\rm Dist}(B_{r,\mathbb{F}_p}^+),\ \ 
\mathcal{U}_r^{\leq 0}={\rm Dist}(B_{r,\mathbb{F}_p}),\]
where $H_{\mathbb{F}_p}$ denotes the base change to 
$\mathbb{F}_p$ of a group scheme $H$ defined over 
$\mathbb{F}_p$. 
For detailed definition of hyperalgebras 
(which are also called algebras of distributions), see \cite[I, \S 7]{jantzenbook}. 
Then we have triangular decompositions 
$\mathcal{U}=\mathcal{U}^-\mathcal{U}^0\mathcal{U}^+$ and 
$\mathcal{U}_r=\mathcal{U}_r^-\mathcal{U}_r^0\mathcal{U}_r^+$ 
for $r \in \mathbb{Z}_{>0}$. In other words, the multiplication 
maps 
$\mathcal{U}^- \otimes_{\mathbb{F}_p} \mathcal{U}^0 \otimes_{\mathbb{F}_p} 
\mathcal{U}^+ \rightarrow \mathcal{U}$ 
and 
$\mathcal{U}_r^- \otimes_{\mathbb{F}_p} \mathcal{U}_r^0 \otimes_{\mathbb{F}_p} 
\mathcal{U}_r^+ \rightarrow \mathcal{U}_r$ 
are   $\mathbb{F}_p$-linear isomorphisms.

Let ${\rm Fr}:\mathcal{U} \rightarrow \mathcal{U}$ be an $\mathbb{F}_p$-algebra 
endomorphism defined by 
\[
{\rm Fr} \left(e_{\alpha}^{(n)} \right)= 
\left\{ \begin{array}{ll}
{e_{\alpha}^{(n/p)}} & \mbox{if $p\ |\ n$}, \\
0 &  \mbox{if $p \nmid n$} \\
\end{array} \right. 
\]
for $\alpha \in \Phi$, which is called the Frobenius endomorphism. 
On the other hand,  there are 
$\mathbb{F}_p$-algebra endomorphisms 
${\rm Fr}'_+:\mathcal{U}^+ \rightarrow \mathcal{U}^+$,  
${\rm Fr}'_-:\mathcal{U}^- \rightarrow \mathcal{U}^-$, and 
${\rm Fr}'_0:\mathcal{U}^0 \rightarrow \mathcal{U}^0$  defined by 
${\rm Fr}'_+(e_{\alpha_i}^{(n)}) = e_{\alpha_i}^{(np)}$, 
${\rm Fr}'_-(e_{-\alpha_i}^{(n)}) = e_{-\alpha_i}^{(np)}$, and 
${\rm Fr}'_0({h_i \choose n}) = {h_i \choose np}$ for $1 \leq i \leq l$  
 and 
$n \in \mathbb{Z}_{\geq 0}$ (see Proposition 1.1 and Corollaire 1.2 in 
\cite{gros-kaneda11}).   Then there is an $\mathbb{F}_p$-linear map 
${\rm Fr}': \mathcal{U} \rightarrow \mathcal{U}$ defined by 
${\rm Fr}'({\bm f} {\bm h} {\bm e})= 
{\rm Fr}'_-({\bm f}) {\rm Fr}'_0({\bm h}){\rm Fr}'_+({\bm e})$ for 
${\bm f} \in \mathcal{U}^-$, ${\bm h} \in \mathcal{U}^0$, and 
${\bm e} \in \mathcal{U}^+$. 
Clearly we have ${\rm Fr} \circ {\rm Fr}' = {\rm id}_{\mathcal{U}}$. 
Since all the maps ${\rm Fr}_{\varepsilon}'$ with 
$\varepsilon \in \{+, -, 0\}$ are restriction ones of ${\rm Fr}'$, from now on 
we denote them by ${\rm Fr}'$ again.

The following proposition is well-known as Lucas' theorem, which is necessary in 
carrying out calculation in $\mathbb{F}_p$. \\ 
\ 
  
\begin{Prop}\label{binomial1}
Let $m, n \in \mathbb{Z}_{\geq 0}$. Let $m= \sum_{i \geq 0} m_i p^i$ and 
$n= \sum_{i \geq 0} n_i p^i$ be their $p$-adic expansions. Then we have 
\[ {m\choose n} \equiv \prod_{i \geq 0}{m_i \choose n_i}\ ({\rm mod}\ p).\]
\end{Prop}
\ 

\noindent {\bf Remark.} Let $r \in \mathbb{Z}_{>0}$ and $\alpha \in \Phi$. 
If we write a positive integer $n \in \mathbb{Z}_{>0}$ as 
$n= n' + n''p^r$ with $n' \in \mathcal{N}_r$ and $n'' \in \mathbb{Z}_{\geq 0}$, then 
by Proposition \ref{basicformulas} (v) and Proposition \ref{binomial1} we have 
\[
{h_{\alpha} \choose n} ={h_{\alpha} \choose n'} {h_{\alpha} \choose n''p^r}
\] 
in $\mathcal{U}$. 
\\ \\ \  

For later use, we also give the following fact. \\ 
\ 

\begin{Prop}\label{binomialforh}
Let $\alpha \in \Phi$, $c,m,n \in \mathbb{Z}$, and 
$r \in \mathbb{Z}_{>0}$. If  
$0 \leq n \leq p^r-1$, then 
\[
{h_{\alpha}+c+mp^r \choose n} = {h_{\alpha}+c \choose n}
\]
in $\mathcal{U}$. 
\end{Prop}

\noindent {\itshape Proof.} 
Since $c$ is arbitrarily chosen, we may assume that $m \geq 0$ without loss of 
generality. 
It is well-known that 
\[
{h_{\alpha}+c+mp^r \choose n} = \sum_{i=0}^{n} 
{mp^r \choose i} {h_{\alpha}+c \choose n-i}
\]
in $\mathcal{U}_{\mathbb{Z}}$, and hence in $\mathcal{U}$ 
(see \cite[Corollary 3.1.2]{gros12}). 
Then by Proposition \ref{binomial1} and the assumption on $n$ 
we must have 
\[
{mp^r \choose i} =0
\]
in $\mathbb{F}_p$ for $0 < i \leq n$. Thus the result follows. 
$\square$ \\

Let $K$ be a field and $R_0, R_1, \dots$ a sequence of finite-dimensional 
associative $K$-algebras with $1$. For each $i \in \mathbb{Z}_{\geq 0}$, let $V_i$ be 
a $K$-subspace of $R_i$ containing $1$. For $n \in \mathbb{Z}_{\geq 0}$, consider 
the subspace $V(n)=V_0 \otimes_{K} V_1 \otimes_K \cdots \otimes_K V_n$ of the 
$K$-algebra $R(n)=R_0 \otimes_{K} R_1 \otimes_K \cdots \otimes_K R_n$. If we 
regard $R(n)$ as a $K$-subalgebra of $R(n+1)$ via 
\[
R(n) \rightarrow R(n+1), \ \ x_0 \otimes \cdots \otimes x_{n} \mapsto 
x_0 \otimes \cdots \otimes x_{n} \otimes 1,
\]
where $x_i \in R_i$, then 
we can regard $V(n)$ as a $K$-subspace of $V(n+1)$. Now we define a 
$K$-vector space $\bigotimes_{n \geq 0} V_n$ as the union 
$\bigcup_{n \geq 0} V(n)$.

\section{Commutation formulas in $\mathcal{U}_{\mathbb{Z}}$ and $\mathcal{U}$}

In this section, we give various commutation formulas in 
$\mathcal{U}_{\mathbb{Z}}$ and $\mathcal{U}$, which will be used to prove main results. 
We first consider such formulas on  $e_{\alpha}$ and $e_{\beta}$ for two roots 
$\alpha, \beta \in \Phi$ with $\alpha + \beta \in \Phi$. 
We know that the subset 
$\Phi'(\alpha, \beta)=(\mathbb{Z} \alpha+\mathbb{Z} \beta) \cap \Phi$ forms a root  
system of type ${\rm A}_2$, ${\rm B}_2$, 
or ${\rm G}_2$. 
If we choose  a (unique) 
nonnegative integer $m$ such that $\beta-m \alpha \in \Phi$ and 
$\beta-(m+1) \alpha \not\in \Phi$, then there exists  $c_{\alpha, \beta} \in \{ \pm 1\}$ 
such that $[e_{\alpha}, e_{\beta}]=(m+1) c_{\alpha, \beta} e_{\alpha+\beta}$ in 
$\mathfrak{g}_{\mathbb{Z}}$ (see \cite[Theorem 25.2 (d)]{humphreysbook}).  
From now on, we often use the symbol $c_{\alpha, \beta}$ to express such a number. 
Without loss of generality, we assume for a moment that $||\alpha|| \leq ||\beta||$ and 
$\alpha$ and $\beta$ form a base of the subsystem $\Phi'(\alpha, \beta)$. 

Suppose that $\Phi'(\alpha, \beta)$ is of 
type ${\rm A}_2$. Then we have $||\beta||=||\alpha||$, 
\[
\Phi'(\alpha, \beta)=\{ \pm \alpha, \pm \beta, \pm(\alpha+\beta)\}, 
\]
and 
\[
\angle(\alpha, \beta)=2\pi/3,\ \   
\angle(\alpha, \alpha+\beta)=\pi/3.
\] 
We can write 
$[e_{\alpha}, e_{\beta}]= c_{\alpha, \beta} e_{\alpha+\beta}$ in 
$\mathfrak{g}_{\mathbb{Z}}$ for some $c_{\alpha, \beta} \in \{ \pm 1\}$. Then for 
$a,b \in \mathbb{Z}_{\geq 0}$, using induction we see that 
\begin{align}\label{commforma1} e_{\alpha}^{(a)} e_{\beta}^{(b)} = 
\sum_{\substack{t_1+t_2=b, \\ t_2+t_3=a}} c_{\alpha, \beta}^{t_2} 
e_{\beta}^{(t_1)} e_{\alpha+\beta}^{(t_2)} 
e_{\alpha}^{(t_3)},
\end{align}
\begin{align}\label{commforma2}  e_{\beta}^{(b)} e_{\alpha}^{(a)} = 
\sum_{\substack{t_1+t_2=a, \\ t_2+t_3=b}} (-c_{\alpha, \beta})^{t_2} 
e_{\alpha}^{(t_1)} e_{\alpha+\beta}^{(t_2)} 
e_{\beta}^{(t_3)}
\end{align} 
in $\mathcal{U}_{\mathbb{Z}}$. 

Suppose that $\Phi'(\alpha, \beta)$ is of 
type ${\rm B}_2$. Then 
$||\beta||=\sqrt{2}||\alpha||$, 
\[
\Phi'(\alpha, \beta)=\{ \pm \alpha, \pm \beta, \pm(\alpha+\beta), \pm (2\alpha+\beta)\}, 
\] 
and 
\[ \angle(\alpha, \beta)=3\pi/4,\ \  
\angle(\alpha, \alpha+\beta)=\pi/2, \ \   \angle(\alpha, 2\alpha+\beta)=\pi/4.\]  
We can write 
$[e_{\alpha}, e_{\beta}]= c_{\alpha,\beta} e_{\alpha+\beta}$ and 
$[e_{\alpha}, e_{\alpha+\beta}] = 2c_{\alpha, \alpha+\beta} e_{2\alpha+\beta}$ in 
$\mathfrak{g}_{\mathbb{Z}}$ for some 
$c_{\alpha,\beta}, c_{\alpha, \alpha+\beta} \in \{ \pm 1\}$. Then for 
$a,b \in \mathbb{Z}_{\geq 0}$, using induction we see that 
\begin{align}\label{commformb1}  e_{\alpha}^{(a)} e_{\beta}^{(b)} = 
\sum_{\substack{t_1+t_2+t_3=b, \\ t_2+2t_3+t_4=a}} c_{\alpha,\beta}^{t_2} 
(c_{\alpha, \beta}c_{\alpha, \alpha+\beta})^{t_3} 
e_{\beta}^{(t_1)} e_{\alpha+\beta}^{(t_2)} e_{2\alpha+\beta}^{(t_3)} e_{\alpha}^{(t_4)},
\end{align}
\begin{align}\label{commformb2}  e_{\beta}^{(b)} e_{\alpha}^{(a)} = 
\sum_{\substack{t_1+2t_2+t_3=a, \\ t_2+t_3+t_4=b}} 
(-c_{\alpha, \beta})^{t_3} (c_{\alpha, \beta} c_{\alpha, \alpha+\beta})^{t_2} 
e_{\alpha}^{(t_1)} e_{2\alpha+\beta}^{(t_2)} e_{\alpha+\beta}^{(t_3)} e_{\beta}^{(t_4)},
\end{align} 
\begin{align}\label{commformb3}  e_{\alpha}^{(a)} e_{\alpha+\beta}^{(b)} = 
\sum_{\substack{t_1+t_2=b, \\ t_2+t_3=a}} (2c_{\alpha, \alpha+\beta})^{t_2} 
e_{\alpha+\beta}^{(t_1)} e_{2\alpha+\beta}^{(t_2)} e_{\alpha}^{(t_3)},
\end{align} 
\begin{align}\label{commformb4}  e_{\alpha+\beta}^{(b)} e_{\alpha}^{(a)} = 
\sum_{\substack{t_1+t_2=a, \\ t_2+t_3=b}} (-2c_{\alpha, \alpha+\beta})^{t_2} 
e_{\alpha}^{(t_1)} e_{2\alpha+\beta}^{(t_2)} e_{\alpha+\beta}^{(t_3)}
\end{align} 
in $\mathcal{U}_{\mathbb{Z}}$. 

Suppose that $\Phi'(\alpha, \beta)$ is 
of type ${\rm G}_2$. Then 
$||\beta||=\sqrt{3}||\alpha||$, 
\[
\Phi'(\alpha, \beta)=\{ \pm \alpha, \pm \beta, \pm(\alpha+\beta), \pm (2\alpha+\beta), 
\pm (3\alpha+\beta), \pm (3\alpha+2\beta)\}, 
\] 
and  
\[ \angle(\alpha, \beta)=5\pi/6,\ \  
\angle(\alpha, \alpha+\beta)=\angle(3\alpha+\beta, \beta)=2\pi/3, \]
\[ \angle(\alpha, 2\alpha+\beta)=
\angle(2\alpha+\beta, \alpha+\beta)=\pi/3.\] 
We can write 
$[e_{\alpha}, e_{\beta}]= c_{\alpha,\beta} e_{\alpha+\beta}$,  
$[e_{\alpha}, e_{\alpha+\beta}] = 2c_{\alpha, \alpha+\beta} e_{2\alpha+\beta}$, 
$[e_{\alpha}, e_{2\alpha+\beta}] = 3c_{\alpha, 2\alpha+\beta} e_{3\alpha+\beta}$, 
$[e_{2\alpha+\beta}, e_{\alpha+\beta}]
= 3c_{2\alpha+\beta,\alpha+\beta} e_{3\alpha+2\beta}$ 
in $\mathfrak{g}_{\mathbb{Z}}$ for some 
\[
c_{\alpha,\beta}, c_{\alpha, \alpha+\beta}, c_{\alpha, 2\alpha+\beta}, 
c_{2\alpha+\beta,\alpha+\beta} \in \{ \pm 1\}. 
\]
Then we also 
have $[e_{3\alpha+\beta}, e_{\beta}] = -c_{\alpha, \beta} c_{\alpha, 2\alpha+\beta} 
c_{2\alpha+\beta,\alpha+\beta} e_{3\alpha+2\beta}$ 
in $\mathfrak{g}_{\mathbb{Z}}$. 
For 
$a,b \in \mathbb{Z}_{\geq 0}$, using induction we see that 
\begin{align}\label{commformg1}  
e_{\alpha}^{(a)} e_{\beta}^{(b)}
= 
\sum_{\substack{t_1+t_2+2t_3+t_4+t_5=b, \\ t_2+3t_3+2t_4+3t_5+t_6=a}} 
 d_1(t_2,t_3,t_4,t_5)
e_{\beta}^{(t_1)} e_{\alpha+\beta}^{(t_2)} e_{3\alpha+2\beta}^{(t_3)} e_{2\alpha+\beta}^{(t_4)}
e_{3\alpha+\beta}^{(t_5)} e_{\alpha}^{(t_6)},
\end{align}
\begin{align}\label{commformg2}  
e_{\beta}^{(b)} e_{\alpha}^{(a)}= 
\sum_{\substack{t_1+3t_2+2t_3+3t_4+t_5=a, \\ t_2+t_3+2t_4+t_5+t_6=b}} 
d_2(t_2,t_3,t_4,t_5)
e_{\alpha}^{(t_1)} e_{3\alpha+\beta}^{(t_2)} e_{2\alpha+\beta}^{(t_3)} 
e_{3\alpha+2\beta}^{(t_4)} e_{\alpha+\beta}^{(t_5)} e_{\beta}^{(t_6)},
\end{align} 
\begin{align}\label{commformg3}  e_{\alpha}^{(a)} e_{\alpha+\beta}^{(b)} = 
\sum_{\substack{t_1+2t_2+t_3+t_4=b, \\ t_2+t_3+2t_4+t_5=a}} d_3(t_2,t_3,t_4)
e_{\alpha+\beta}^{(t_1)} e_{3\alpha+2\beta}^{(t_2)} e_{2\alpha+\beta}^{(t_3)} 
e_{3\alpha+\beta}^{(t_4)}e_{\alpha}^{(t_5)}, 
\end{align} 
\begin{align}\label{commformg4}  e_{\alpha+\beta}^{(b)} e_{\alpha}^{(a)} = 
\sum_{\substack{t_1+2t_2+t_3+t_4=a, \\ t_2+t_3+2t_4+t_5=b}} d_4(t_2,t_3,t_4)
e_{\alpha}^{(t_1)} e_{3\alpha+\beta}^{(t_2)} e_{2\alpha+\beta}^{(t_3)} 
e_{3\alpha+2\beta}^{(t_4)}e_{\alpha+\beta}^{(t_5)}, 
\end{align} 
\begin{align}\label{commformg5} e_{\alpha}^{(a)} e_{2\alpha+\beta}^{(b)} = 
\sum_{\substack{t_1+t_2=b, \\ t_2+t_3=a}} 
(3c_{\alpha,2\alpha+\beta})^{t_2} e_{2\alpha+\beta}^{(t_1)} e_{3\alpha+\beta}^{(t_2)} 
e_{\alpha}^{(t_3)},
\end{align}
\begin{align}\label{commformg6} e_{2\alpha+\beta}^{(b)} e_{\alpha}^{(a)}  = 
\sum_{\substack{t_1+t_2=a, \\ t_2+t_3=b}} (-3c_{\alpha,2\alpha+\beta})^{t_2} e_{\alpha}^{(t_1)} e_{3\alpha+\beta}^{(t_2)} 
e_{2\alpha+\beta}^{(t_3)},
\end{align}
\begin{align}\label{commformg7}  e_{2\alpha+\beta}^{(a)} e_{\alpha+\beta}^{(b)} = 
\sum_{\substack{t_1+t_2=b, \\ t_2+t_3=a}} (3c_{2\alpha+\beta, \alpha+\beta})^{t_2} 
e_{\alpha+\beta}^{(t_1)} e_{3\alpha+2\beta}^{(t_2)} e_{2\alpha+\beta}^{(t_3)}
\end{align} 
\begin{align}\label{commformg8} e_{\alpha+\beta}^{(b)}  e_{2\alpha+\beta}^{(a)} = 
\sum_{\substack{t_1+t_2=a, \\ t_2+t_3=b}} (-3c_{2\alpha+\beta, \alpha+\beta})^{t_2} 
e_{2\alpha+\beta}^{(t_1)} e_{3\alpha+2\beta}^{(t_2)} e_{\alpha+\beta}^{(t_3)}
\end{align} 
\begin{align}\label{commformg9}  e_{3\alpha+\beta}^{(a)} e_{\beta}^{(b)} = 
\sum_{\substack{t_1+t_2=b, \\ t_2+t_3=a}} 
(-c_{\alpha,\beta} c_{\alpha,2\alpha+\beta} c_{2\alpha+\beta, \alpha+\beta})^{t_2} 
e_{\beta}^{(t_1)} e_{3\alpha+2\beta}^{(t_2)} e_{3\alpha+\beta}^{(t_3)}
\end{align} 
\begin{align}\label{commformg10} e_{\beta}^{(b)}  e_{3\alpha+\beta}^{(a)} = 
\sum_{\substack{t_1+t_2=a, \\ t_2+t_3=b}} 
(c_{\alpha,\beta} c_{\alpha,2\alpha+\beta} c_{2\alpha+\beta, \alpha+\beta})^{t_2} 
e_{3\alpha+\beta}^{(t_1)} e_{3\alpha+2\beta}^{(t_2)} e_{\beta}^{(t_3)}
\end{align}
in $\mathcal{U}_{\mathbb{Z}}$, where 
\[ d_1(t_2,t_3,t_4,t_5) = c_{\alpha,\beta}^{t_2} 
(c_{\alpha,\beta}c_{\alpha,\alpha+\beta})^{t_4} 
(c_{\alpha,\beta}c_{\alpha,\alpha+\beta} c_{\alpha,2\alpha+\beta})^{t_5} 
(c_{\alpha,\alpha+\beta}c_{2\alpha+\beta,\alpha+\beta})^{t_3}, \] 
\[ d_2(t_2,t_3,t_4,t_5) =(-c_{\alpha,\beta})^{t_5} (c_{\alpha,\beta}c_{\alpha,\alpha+\beta})^{t_3} 
(-c_{\alpha,\beta}c_{\alpha,\alpha+\beta} c_{\alpha,2\alpha+\beta})^{t_2} 
(c_{\alpha,\alpha+\beta}c_{2\alpha+\beta,\alpha+\beta})^{t_4}, \] 
\[ d_3(t_2,t_3,t_4) =(2c_{\alpha,\alpha+\beta})^{t_3} 
(3c_{\alpha,\alpha+\beta}c_{\alpha,2\alpha+\beta})^{t_4} 
(3c_{\alpha,\alpha+\beta}c_{2\alpha+\beta,\alpha+\beta})^{t_2}, \] 
\[ d_4(t_2,t_3,t_4) =(-2c_{\alpha,\alpha+\beta})^{t_3} 
(3c_{\alpha,\alpha+\beta}c_{\alpha,2\alpha+\beta})^{t_2} 
(3c_{\alpha,\alpha+\beta}c_{2\alpha+\beta,\alpha+\beta})^{t_4}.\] 

Note that the sums in the right-hand sides of the equalities 
(\ref{commforma1})-(\ref{commformg10}) are finite (recall that $e_{\gamma}^{(n)}=0$ 
for $\gamma \in \Phi$ and $n <0$). \\

Now let $\alpha, \gamma \in \Phi$ be two roots such that $\alpha + \gamma \in \Phi$. 
But we do not assume that $\alpha$ and $\gamma$ form a base of 
 $\Phi'(\alpha, \gamma) = (\mathbb{Z}\alpha + \mathbb{Z}\gamma) \cap \Phi$.  
The possible cases are as follows:
\[ \begin{array}{|c||c|c|c|c|} \hline 
& \angle(\alpha, \gamma) & \Phi'(\alpha, \gamma) 
 & \mbox{lengths of $\alpha$, $\gamma$, $\alpha+\gamma$} & m \\ \hline \hline 
{\rm (A)} & 2\pi/3   & {\rm A}_2 & ||\alpha+\gamma||=||\alpha||=||\gamma|| &  0 \\ \hline 
{\rm (B)} & \pi/2   & {\rm B}_2 & ||\alpha+\gamma||=\sqrt{2}||\alpha||=
\sqrt{2}||\gamma|| & 1 \\ \hline 
{\rm (C)} & \pi/3   & {\rm G}_2 & ||\alpha+\gamma||=\sqrt{3}||\alpha||
=\sqrt{3}||\gamma|| & 2 \\ \hline 
{\rm (D)} & 3\pi/4   & {\rm B}_2 & \sqrt{2}||\alpha+\gamma||
=\sqrt{2}||\alpha||=||\gamma|| & 0 \\ \hline
{\rm (E)} & 3\pi/4   & {\rm B}_2 & \sqrt{2}||\alpha+\gamma||
=||\alpha||=\sqrt{2}||\gamma|| & 0 \\ \hline
{\rm (F)} & 5\pi/6   & {\rm G}_2 & \sqrt{3}||\alpha+\gamma||
=\sqrt{3}||\alpha||=||\gamma|| & 0 \\ \hline 
{\rm (G)} & 5\pi/6   & {\rm G}_2 & \sqrt{3}||\alpha+\gamma||
=||\alpha||=\sqrt{3}||\gamma|| & 0 \\ \hline 
{\rm (H)} & 2\pi/3   & {\rm G}_2 & ||\alpha+\gamma||
=||\alpha||=||\gamma|| 
& 1 \\ \hline
\end{array}\]
\centerline{\bf Table 1} \\ 
Here $m$ is a unique nonnegative integer satisfying $\gamma - m\alpha \in \Phi$ and 
$\gamma - (m+1)\alpha \not\in \Phi$. 

The above table will be used later. \\

\begin{Lem}\label{binomial2}
Let $n, r\in \mathbb{Z}_{> 0}$ and $m \in \mathbb{Z}$ with  $0 \leq m \leq np^r$. 
Then in $\mathbb{F}_p$ we have
\[ \sum_{i=0}^n (-1)^i {np^r-m \choose ip^r} =
\left\{ \begin{array}{ll}
1 & \mbox{if $(n-1)p^r < m \leq np^r$,} \\
0 & \mbox{if $m \leq (n-1)p^r$}
\end{array} \right.. \] 
\end{Lem}

\noindent {\itshape Proof.} If $(n-1)p^r < m \leq np^r$, then 
$0 \leq np^r-m< p^r$ and hence 
we have 
\[
\sum_{i=0}^n (-1)^i {np^r-m \choose ip^r} =(-1)^0 {np^r-m \choose 0}=1
\]
in $\mathbb{Z}$, and hence in $\mathbb{F}_p$. 
So assume that $m \leq (n-1)p^r$. Then $np^r-m=n' p^r+n''$ for some integers $n'$ and 
$n''$ with $0 < n' \leq n$ and $n'' \in \mathcal{N}_r$. By Proposition \ref{binomial1}, we have 
\[
\sum_{i=0}^n (-1)^i {np^r-m \choose ip^r} =\sum_{i=0}^n (-1)^i {n'p^r+n'' \choose ip^r} =
\sum_{i=0}^{n'} (-1)^i {n' \choose i} =(1-1)^{n'}=0 
\]
 in $\mathbb{F}_p$, as required. $\square$ \\
\ 

\begin{Prop}\label{commformpr1}
Let $\alpha, \gamma \in \Phi$ be two roots such that $\alpha + \gamma \in \Phi$. 
Then for $n, r \in \mathbb{Z}_{>0}$, the element 
$z=\sum_{i=0}^n (-1)^i e_{\alpha}^{((n-i)p^r)} e_{\gamma}^{(p^r)} e_{\alpha}^{(ip^r)}$ in 
$\mathcal{U}$ is as follows:
\renewcommand{\multirowsetup}{\centering}
\hfill
\[ \begin{array}{|c||c|c|} \hline 
\mbox{Case in Table 1} & n & z
  \\ \hline \hline 
 \multirow{2}{30mm}{\mbox{{\rm (A)}, {\rm (E)}, or {\rm (G)}}} & n \geq 2   & 0   \\ 
  & n=1 & c_{\alpha, \gamma} e_{\alpha+\gamma}^{(p^r)}+
(\mbox{an element of $\mathcal{U}_r$}) 
\\ \hline
 \multirow{2}{30mm}{\rm (B)} & n \geq 2   & 0   \\ 
  & n=1 & 2c_{\alpha, \gamma} e_{\alpha+\gamma}^{(p^r)}
+(\mbox{an element of $\mathcal{U}_r$}) 
\\ \hline
\multirow{2}{30mm}{\rm (C)} & n \geq 2   & 0   \\
  & n=1 & 3c_{\alpha, \gamma} e_{\alpha+\gamma}^{(p^r)}
+(\mbox{an element of $\mathcal{U}_r$}) 
\\ \hline
\multirow{3}{30mm}{\rm (D)} & n \geq 3   & 0   \\ 
   & n=2 & c_{\alpha, \gamma}  c_{\alpha, \alpha+\gamma} 
e_{2\alpha+\gamma}^{(p^r)}+(\mbox{an element of $\mathcal{U}_r$}) \\ 
  & n=1 & c_{\alpha, \gamma} e_{\alpha+\gamma}^{(p^r)}
+(\mbox{an element of $\mathcal{U}_r$}) 
\\ \hline
\multirow{4}{30mm}{\rm (F)} & n \geq 4   & 0   \\ 
   & n=3 & c_{\alpha, \gamma}  c_{\alpha, \alpha+\gamma} c_{\alpha, 2\alpha+\gamma} 
e_{3\alpha+\gamma}^{(p^r)}+(\mbox{an element of $\mathcal{U}_r$}) \\ 
  & n=2 & c_{\alpha, \gamma}  c_{\alpha, \alpha+\gamma} e_{2\alpha+\gamma}^{(p^r)}
+(\mbox{an element of $\mathcal{U}_r$}) \\ 
 & n=1 & c_{\alpha, \gamma}  e_{\alpha+\gamma}^{(p^r)}
+(\mbox{an element of $\mathcal{U}_r$}) 
\\ \hline
\multirow{3}{30mm}{\rm (H)} & n \geq 3   & 0   \\ 
   & n=2 & 3c_{\alpha, \gamma}  c_{\alpha, \alpha+\gamma} 
e_{2\alpha+\gamma}^{(p^r)}+(\mbox{an element of $\mathcal{U}_r$}) \\ 
  & n=1 & 2c_{\alpha, \gamma} e_{\alpha+\gamma}^{(p^r)}
+(\mbox{an element of $\mathcal{U}_r$}) 
\\ \hline
\end{array} \] 
\centerline{\bf Table 2}
In particular, if $n=1$, then we have 
\[
e_{\alpha}^{(p^r)}e_{\gamma}^{(p^r)} - e_{\gamma}^{(p^r)}e_{\alpha}^{(p^r)}= 
k c_{\alpha, \gamma} e_{\alpha+\gamma}^{(p^r)}+(\mbox{an element of $\mathcal{U}_r$})
\]
in $\mathcal{U}$, where  $k \in \{  1,  2,  3\}$ satisfies 
$[e_{\alpha}, e_{\gamma}] = k c_{\alpha, \gamma} e_{\alpha+\gamma}$. 
\end{Prop}
\ 

\noindent {\itshape Proof.} 
The second statement follows from Table 2 and the fact that 
\[
[e_{\alpha}, e_{\gamma}] = (m+1) c_{\alpha, \gamma}e_{\alpha+\gamma}
\]
for $m$ in Table 1. 
The results in Table 2 follow from direct calculation using the 
commutation formulas given before. In all cases, the arguments are similar. We left 
almost all of the proof to the reader, but we shall deal with only the case (F) as an example  
(the other cases are easier). 

Consider the case (F). Note that 
\[
e_{\alpha} e_{\gamma} -e_{\gamma}e_{\alpha} =c_{\alpha, \gamma} e_{\alpha+\gamma},
\]
\[
e_{\alpha} e_{\alpha+\gamma} -e_{\alpha+\gamma}e_{\alpha} =
2c_{\alpha, \alpha+\gamma} e_{2\alpha+\gamma},
\]
\[
e_{\alpha} e_{2\alpha+\gamma} -e_{2\alpha+\gamma}e_{\alpha} =
3c_{\alpha, 2\alpha+\gamma} e_{3\alpha+\gamma},
\]
\[
e_{2\alpha+\gamma} e_{\alpha+\gamma} -e_{\alpha+\gamma}e_{2\alpha+\gamma} =
3c_{2\alpha+\gamma, \alpha+\gamma} e_{3\alpha+2\gamma}.
\]
Then using the formula (\ref{commformg1}) and Proposition \ref{basicformulas} (i) we have 
\begin{align*}
z & = \sum_{i=0}^{n} (-1)^i e_{\alpha}^{((n-i)p^r)} e_{\gamma}^{(p^r)} e_{\alpha}^{(ip^r)} \\
&= \sum_{i=0}^{n} (-1)^i 
\sum_{\substack{t_1+t_2+2t_3+t_4+t_5=p^r, \\ t_2+3t_3+2t_4+3t_5+t_6=(n-i)p^r}} 
 d(t_2,t_3,t_4,t_5)
e_{\gamma}^{(t_1)} e_{\alpha+\gamma}^{(t_2)} e_{3\alpha+2\gamma}^{(t_3)}
 e_{2\alpha+\gamma}^{(t_4)}e_{3\alpha+\gamma}^{(t_5)} e_{\alpha}^{(t_6)} e_{\alpha}^{(ip^r)} \\
&=\sum_{t_1+t_2+2t_3+t_4+t_5=p^r} d(t_2,t_3,t_4,t_5) 
e_{\gamma}^{(t_1)} e_{\alpha+\gamma}^{(t_2)} e_{3\alpha+2\gamma}^{(t_3)}
 e_{2\alpha+\gamma}^{(t_4)}e_{3\alpha+\gamma}^{(t_5)} \\
& \ \ \ \times \sum_{i=0}^n (-1)^i {np^r-t_2-3t_3-2t_4-3t_5 \choose ip^r} 
e_{\alpha}^{(np^r-t_2-3t_3-2t_4-3t_5)},
\end{align*}
in $\mathcal{U}$, where 
\[ d(t_2,t_3,t_4,t_5) = 
c_{\alpha,\gamma}^{t_2} (c_{\alpha,\gamma}c_{\alpha,\alpha+\gamma})^{t_4} 
(c_{\alpha,\gamma}c_{\alpha,\alpha+\gamma}c_{\alpha,2\alpha+\gamma})^{t_5}
(c_{\alpha,\alpha+\gamma} c_{2\alpha+\gamma,\alpha+\gamma})^{t_3}. \]
By Lemma \ref{binomial2},  we have 
\[
\sum_{i=0}^n (-1)^i {np^r-t_2-3t_3-2t_4-3t_5 \choose ip^r} =
\left\{ \begin{array}{ll}
1 & \mbox{if $(n-1)p^r < t_2+3t_3+2t_4+3t_5 \leq np^r$,} \\
0 & \mbox{if $t_2+3t_3+2t_4+3t_5 \leq (n-1)p^r$} 
\end{array} \right.
\] 
in $\mathbb{F}_p$. Note that  $t_1+t_2+2t_3+t_4+t_5=p^r$ implies 
$t_2+3t_3+2t_4+3t_5 \leq 3p^r$. Thus if $n \geq 4$, we must have $z=0$. 
Suppose that $n=1$, $2$, or $3$. In each case, there exists a unique $5$-tuple 
$(t_1,t_2,t_3,t_4,t_5)$ such that one of the $t_i$'s is equal to $p^r$ when 
\[(n-1)p^r < t_2+3t_3+2t_4+3t_5 \leq np^r.\] 
Indeed, if $n=3$, then the $5$-tuple is  $(t_1,t_2,t_3,t_4,t_5)=(0,0,0,0,p^r)$ and the corresponding term in $\mathcal{U}$ 
is $c_{\alpha,\gamma}c_{\alpha,\alpha+\gamma}c_{\alpha,2\alpha+\gamma}
e_{3\alpha+\gamma}^{(p^r)}$. If $n=2$, then the $5$-tuple is  $(t_1,t_2,t_3,t_4,t_5)=(0,0,0,p^r,0)$ and the corresponding term in $\mathcal{U}$ 
is $c_{\alpha,\gamma}c_{\alpha,\alpha+\gamma}e_{2\alpha+\gamma}^{(p^r)}$. 
If $n=1$, then the $5$-tuple is  $(t_1,t_2,t_3,t_4,t_5)=(0,p^r,0,0,0)$ and the corresponding term  in $\mathcal{U}$ 
is $c_{\alpha,\gamma}e_{\alpha+\gamma}^{(p^r)}$. 
On the other hand, all the terms for other $5$-tuples $(t_1,t_2,t_3,t_4,t_5)$ 
lie in $\mathcal{U}_r$ since in this case 
each $t_i$ is less than $p^r$. Now the proposition for the case (F) is proved. $\square$ \\

\begin{Prop}\label{commformpr2}
Let $\alpha \in \Phi$, $n \in \mathbb{Z}_{\geq 0}$, 
$r \in \mathbb{Z}_{> 0}$, and $z \in \mathcal{U}_r$. Then the  element 
\[ \sum_{i=0}^{n} (-1)^i e_{\alpha}^{((n-i)p^r)} z e_{\alpha}^{(ip^r)}\] 
of $\mathcal{U}$ lies in $\mathcal{U}_r$. 
\end{Prop}
\ 

\noindent {\bf Remark.} More generally, for $\alpha \in \Phi$, 
$n, c \in \mathbb{Z}_{\geq 0}$, 
$r \in \mathbb{Z}_{> 0}$, and $z \in \mathcal{U}_r$, 
we see that the element 
\[ \sum_{i=0}^{n} (-1)^i e_{\alpha}^{((n-i)p^r-c)} z e_{\alpha}^{(ip^r)}\] 
of $\mathcal{U}$ lies in $\mathcal{U}_r$. Indeed, 
we may assume that $np^r \geq c$ and then we can write $c=c'p^r+c''$ with 
$0 \leq c' \leq n$ and $c'' \in \mathcal{N}_r$. If $c''=0$, by Proposition 
\ref{commformpr2} we have 
\begin{align*}
\sum_{i=0}^{n} (-1)^i e_{\alpha}^{((n-i)p^r-c)} z e_{\alpha}^{(ip^r)}
&= \sum_{i=0}^{n} (-1)^i e_{\alpha}^{((n-c'-i)p^r)} z e_{\alpha}^{(ip^r)} \\
&= \sum_{i=0}^{n-c'} (-1)^i e_{\alpha}^{((n-c'-i)p^r)} z e_{\alpha}^{(ip^r)} \in \mathcal{U}_r
\end{align*}
in $\mathcal{U}$.  
On the other hand, if $c'' \neq 0$, then $0 \leq c' <n$ and 
by Proposition \ref{basicformulas} (i) and Proposition \ref{commformpr2} we have 
\begin{align*}
\sum_{i=0}^{n} (-1)^i e_{\alpha}^{((n-i)p^r-c)} z e_{\alpha}^{(ip^r)}
&= \sum_{i=0}^{n} (-1)^i e_{\alpha}^{(p^r-c'')}e_{\alpha}^{((n-c'-1-i)p^r)} z e_{\alpha}^{(ip^r)} \\
&= e_{\alpha}^{(p^r-c'')}\sum_{i=0}^{n-c'-1} (-1)^i e_{\alpha}^{((n-c'-1-i)p^r)} z e_{\alpha}^{(ip^r)} \in \mathcal{U}_r
\end{align*}
in $\mathcal{U}$, as required. \\ \

\noindent {\itshape Proof of Proposition \ref{commformpr2}.} 
It is clear when $n=0$, so we may assume that $n >0$. 
We may also assume that $z=\prod_{j=1}^k e_{\gamma_j}^{(m_j)} \in \mathcal{U}_r$ 
with $k \geq 0$, $\gamma_j \in \Phi$, and $m_j \in \mathcal{N}_r-\{0\}$ for each $j$. 
We proceed by induction on $k$. 

Suppose that $k=0$. Then since $z=1$, by 
Proposition \ref{basicformulas} (i) and Lemma \ref{binomial2} we have 
\[
\sum_{i=0}^{n} (-1)^i e_{\alpha}^{((n-i)p^r)} z e_{\alpha}^{(ip^r)}
= \sum_{i=0}^{n} (-1)^i {np^r \choose ip^r} e_{\alpha}^{(np^r)} 
=0 
\]  
in $\mathcal{U}$. 

From now on, we assume that $k >0$. Set 
$z'= \prod_{j=2}^k e_{\gamma_j}^{(m_j)}$, $\gamma=\gamma_1$, and 
$m=m_1$. Then $z=e_{\gamma}^{(m)}z'$. If $\alpha+\gamma \not\in \Phi$ and 
$\gamma \neq -\alpha$, in $\mathcal{U}$ we have 
\[\sum_{i=0}^{n} (-1)^i e_{\alpha}^{((n-i)p^r)} z e_{\alpha}^{(ip^r)}
= e_{\gamma}^{(m)} \sum_{i=0}^{n} (-1)^i e_{\alpha}^{((n-i)p^r)} z' e_{\alpha}^{(ip^r)} 
\in  \mathcal{U}_r\]
by Proposition 2.1 (iv) and induction. If $\gamma = -\alpha$, in $\mathcal{U}$ we have 
\begin{align*}
\lefteqn{\sum_{i=0}^{n} (-1)^i e_{\alpha}^{((n-i)p^r)} z e_{\alpha}^{(ip^r)} = 
 \sum_{i=0}^{n-1} (-1)^i e_{\alpha}^{((n-i)p^r)} e_{-\alpha}^{(m)} z' e_{\alpha}^{(ip^r)}  
+ (-1)^n e_{-\alpha}^{(m)} z' e_{\alpha}^{(np^r)} } \\
&= \sum_{i=0}^{n-1} (-1)^i \sum_{c=0}^m e_{-\alpha}^{(m-c)} 
{h_{\alpha} -m -(n-i)p^r+2c \choose c} e_{\alpha}^{((n-i)p^r-c)} z' e_{\alpha}^{(ip^r)} 
+ (-1)^n e_{-\alpha}^{(m)} z' e_{\alpha}^{(np^r)}  \\
&=  \sum_{i=0}^{n-1} (-1)^i \sum_{c=0}^m e_{-\alpha}^{(m-c)} 
{h_{\alpha} -m +2c \choose c} e_{\alpha}^{((n-i)p^r-c)} z' e_{\alpha}^{(ip^r)} 
+ (-1)^n e_{-\alpha}^{(m)} z' e_{\alpha}^{(np^r)}  \\
&=  \sum_{i=0}^{n} (-1)^i \sum_{c=0}^m e_{-\alpha}^{(m-c)} 
{h_{\alpha} -m +2c \choose c} e_{\alpha}^{((n-i)p^r-c)} z' e_{\alpha}^{(ip^r)} \\
&= \sum_{c=0}^m e_{-\alpha}^{(m-c)} 
{h_{\alpha} -m +2c \choose c} \sum_{i=0}^n (-1)^i 
e_{\alpha}^{((n-i)p^r-c)} z' e_{\alpha}^{(ip^r)} \in \mathcal{U}_r
\end{align*}
by Proposition \ref{basicformulas} (ii), Proposition \ref{binomialforh}, and  induction 
(recall that $e_{\beta}^{(a)} =0$ for $\beta \in \Phi$ and $a<0$). 

From now on,  we assume that $\alpha+ \gamma \in \Phi$. 
Consider first one of the cases (A), (B), and (C) in Table 1. 
Then there is an integer $d \in \{ \pm1, \pm 2,\pm 3\}$ 
such that $[e_{\alpha}, e_{\gamma}]=d e_{\alpha+\gamma}$ in $\mathfrak{g}_{\mathbb{Z}}$ 
and hence 
\[e_{\alpha}^{((n-i)p^r)} e_{\gamma}^{(m)} = 
\sum_{\substack{t_1+t_2=m, \\ t_2+t_3=(n-i)p^r}} 
d^{t_2} e_{\gamma}^{(t_1)} e_{\alpha+\gamma}^{(t_2)} 
e_{\alpha}^{(t_3)}\]
in $\mathcal{U}_{\mathbb{Z}}$. Therefore, in $\mathcal{U}$ we have 
\begin{align*}
\lefteqn{\sum_{i=0}^{n} (-1)^i e_{\alpha}^{((n-i)p^r)} z e_{\alpha}^{(ip^r)} = 
 \sum_{i=0}^{n} (-1)^i e_{\alpha}^{((n-i)p^r)} e_{\gamma}^{(m)} z' e_{\alpha}^{(ip^r)} } \\
&= \sum_{i=0}^n (-1)^i  
\sum_{\substack{t_1+t_2=m, \\ t_2+t_3=(n-i)p^r}} 
d^{t_2} e_{\gamma}^{(t_1)} e_{\alpha+\gamma}^{(t_2)} 
e_{\alpha}^{(t_3)} z' e_{\alpha}^{(ip^r)} \\
&= \sum_{t_1+t_2=m} 
d^{t_2} e_{\gamma}^{(t_1)} e_{\alpha+\gamma}^{(t_2)} \sum_{i=0}^n (-1)^i  
e_{\alpha}^{((n-i)p^r-t_2)} z' e_{\alpha}^{(ip^r)}   \in \mathcal{U}_r
\end{align*}
by  induction.  

In the other cases in Table 1, the arguments are similar. 
So we left the rest of the proof to the reader, 
but we shall only deal with the case (H). 
Then we can write  $[e_{\alpha}, e_{\gamma}]=2c_{\alpha, \gamma} e_{\alpha+\gamma}$, 
$[e_{\alpha}, e_{\alpha+\gamma}]=3c_{\alpha,\alpha+\gamma} e_{2\alpha+\gamma}$, 
and  $[e_{\alpha+\gamma}, e_{\gamma}]
=3c_{\alpha+\gamma,\gamma} e_{\alpha+2\gamma}$ 
in $\mathfrak{g}_{\mathbb{Z}}$ for some 
$c_{\alpha, \gamma}, c_{\alpha, \alpha+\gamma}, 
c_{\alpha+\gamma, \gamma} \in \{ \pm 1\}$. Then  
 in $\mathcal{U}$ we have 
\begin{align*}
\lefteqn{\sum_{i=0}^{n} (-1)^i e_{\alpha}^{((n-i)p^r)} z e_{\alpha}^{(ip^r)} = 
 \sum_{i=0}^{n} (-1)^i e_{\alpha}^{((n-i)p^r)} e_{\gamma}^{(m)} z' e_{\alpha}^{(ip^r)} } \\
&= \sum_{i=0}^n (-1)^i  
\sum_{\substack{t_1+2t_2+t_3+t_4=m, \\ t_2+t_3+2t_4+t_5=(n-i)p^r}} d(t_2,t_3,t_4)
e_{\gamma}^{(t_1)} e_{\alpha+2\gamma}^{(t_2)} e_{\alpha+\gamma}^{(t_3)} 
e_{2\alpha+\gamma}^{(t_4)}e_{\alpha}^{(t_5)} z'  e_{\alpha}^{(ip^r)}  \\
&= \sum_{t_1+2t_2+t_3+t_4=m} d(t_2,t_3,t_4)
e_{\gamma}^{(t_1)} e_{\alpha+2\gamma}^{(t_2)} e_{\alpha+\gamma}^{(t_3)} 
e_{2\alpha+\gamma}^{(t_4)} \sum_{i=0}^n (-1)^i 
e_{\alpha}^{((n-i)p^r-t_2-t_3-2t_4)} z'  e_{\alpha}^{(ip^r)}  \in \mathcal{U}_r
\end{align*}
by equation (\ref{commformg3}) and induction, 
where 
\[
d(t_2,t_3,t_4)= (2c_{\alpha, \gamma})^{t_3} 
(3c_{\alpha, \gamma}c_{\alpha,\alpha+\gamma})^{t_4} 
(3c_{\alpha, \gamma}c_{\alpha+\gamma,\gamma})^{t_2}.
\] 
Thus the result for the case (H) follows.  $\square$ \\
\

Let $\nu$ be the number of elements in $\Phi^+$. 
Let $w_0$ be the longest element of $W$ and fix a reduced expression  
$w_0=s_{i_1} s_{i_2} \cdots s_{i_{\nu}}$, where 
$s_{i}$ denotes the simple reflection $s_{\alpha_{i}}$. If we set 
\[ \beta_1=\alpha_{i_1}, \beta_2=s_{i_1}(\alpha_{i_2}), \dots, 
\beta_{\nu}=s_{i_1} \dots s_{i_{\nu-1}}(\alpha_{i_{\nu}}),\]
then we have $\Phi^+=\{ \beta_1, \beta_2,\dots, \beta_{\nu}\}$ 
(see \cite[5.6 Exercise 1]{humphreysbook2}). 
So the elements 
\[ e_{\beta_1}^{(a_1)} e_{\beta_2}^{(a_2)} \cdots e_{\beta_{\nu}}^{(a_{\nu})}\]
with $a_i \in \mathbb{Z}_{\geq 0}$ for $1 \leq i \leq \nu$ form a $\mathbb{Z}$-basis of 
$\mathcal{U}_{\mathbb{Z}}^+$ and an $\mathbb{F}_p$-basis of $\mathcal{U}^+$ 
(see \cite[II, 1.12]{jantzenbook}). In the rest 
of the paper, we shall keep the above notation.

The following fact is described in \cite[Proposition 3.2]{yoshii22-2}. \\

\begin{Prop}\label{commformbeta}
Suppose that $\nu >1$. 
For $a,b \in \mathbb{Z}_{>0}$ and  
$j,k \in \mathbb{Z}$ with $1 \leq j < k \leq \nu$, the element 
$e_{\beta_k}^{(a)} e_{\beta_j}^{(b)}-e_{\beta_j}^{(b)}e_{\beta_k}^{(a)}$ 
in $\mathcal{U}_{\mathbb{Z}}$ is a $\mathbb{Z}$-linear 
combination of elements of the form $e_{\beta_j}^{(a_j)} \cdots e_{\beta_k}^{(a_k)}$ 
satisfying the following: \\ \\
$\bullet$ $a_j <b$ and $a_k < a$. \\
$\bullet$ $\sum_{i=j}^{k-1} a_i \leq b$ and $\sum_{i=j+1}^{k}a_i \leq a$. \\
\end{Prop}
\ 

\noindent {\bf Remark.} From the first condition on 
$e_{\beta_j}^{(a_j)} \cdots e_{\beta_k}^{(a_k)}$ in the proposition, we also see that the element 
$e_{\beta_{j+1}}^{(a)} e_{\beta_j}^{(b)} - e_{\beta_j}^{(b)} e_{\beta_{j+1}}^{(a)}$ must be zero in 
$\mathcal{U}_{\mathbb{Z}}$ for any $a,b \in \mathbb{Z}_{\geq 0}$ and $j \in \mathbb{Z}$ 
with $1 \leq j \leq \nu-1$. \\ \\
\ 

The fact that $e_{\beta_k}^{(a)} e_{\beta_j}^{(b)}-e_{\beta_j}^{(b)}e_{\beta_k}^{(a)}$ 
in $\mathcal{U}_{\mathbb{Z}}$ is a $\mathbb{Z}$-linear 
combination of elements of the form $e_{\beta_j}^{(a_j)} \cdots e_{\beta_k}^{(a_k)}$ 
in Proposition \ref{commformbeta} implies the following facts. \\

\begin{Prop}\label{subring}
Let $j,k$ be integers satisfying $1 \leq j \leq k \leq \nu$. 
Let $r \in \mathbb{Z}_{> 0}$. Then the following hold. \\

\noindent {\rm (i)} A $\mathbb{Z}$-span of the elements 
$e_{\beta_j}^{(a_j)} \cdots e_{\beta_k}^{(a_k)}$ with 
$(a_j, \dots, a_k) \in (\mathbb{Z}_{\geq 0})^{k-j+1}$ forms a subring of 
$\mathcal{U}_{\mathbb{Z}}^+$. \\

\noindent {\rm (ii)} An $\mathbb{F}_p$-span of the elements 
$e_{\beta_j}^{(a_j)} \cdots e_{\beta_k}^{(a_k)}$ with 
$(a_j, \dots, a_k) \in (\mathbb{Z}_{\geq 0})^{k-j+1}$ forms an $\mathbb{F}_p$-subalgebra of 
$\mathcal{U}^+$. \\

\noindent {\rm (iii)} An $\mathbb{F}_p$-span of the elements 
$e_{\beta_j}^{(a_j)} \cdots e_{\beta_k}^{(a_k)}$ with 
$a_i \in \mathcal{N}_r$ for $j \leq i \leq k$ forms an $\mathbb{F}_p$-subalgebra of 
$\mathcal{U}_r^+$. \\

\noindent {\rm (iv)} Let $e_{\beta_j}^{(a_j)} \cdots e_{\beta_k}^{(a_k)}$ be a fixed element of 
$\mathcal{U}$ satisfying $a_i \in \mathcal{N}_r$ for each $i$ with $j \leq i \leq k$. 
Let $c \in \mathbb{Z}_{> 0}$. Then the following hold. \\  

\noindent $\bullet$ If $k \neq \nu$, then the element 
\[e_{\beta_{k+1}}^{(c)}e_{\beta_j}^{(a_j)} \cdots e_{\beta_k}^{(a_k)} - 
e_{\beta_j}^{(a_j)} \cdots e_{\beta_k}^{(a_k)} e_{\beta_{k+1}}^{(c)}\]
in $\mathcal{U}$ is an $\mathbb{F}_p$-linear combination of elements of the form 
$e_{\beta_j}^{(b_j)} \cdots e_{\beta_k}^{(b_k)} e_{\beta_{k+1}}^{(b_{k+1})}$ satisfying 
$b_{k+1} < c$ and $b_i \in \mathcal{N}_r$ for $j \leq i \leq k$. \\

\noindent $\bullet$ If $j \neq 1$, then the element 
\[e_{\beta_j}^{(a_j)} \cdots e_{\beta_k}^{(a_k)} e_{\beta_{j-1}}^{(c)}- 
e_{\beta_{j-1}}^{(c)}e_{\beta_j}^{(a_j)} \cdots e_{\beta_k}^{(a_k)} \]
in $\mathcal{U}$ is an $\mathbb{F}_p$-linear combination of elements of the form 
$e_{\beta_{j-1}}^{(b_{j-1})}e_{\beta_j}^{(b_j)} \cdots e_{\beta_k}^{(b_k)} $ satisfying 
$b_{j-1} < c$ and $b_i \in \mathcal{N}_r$ for $j \leq i \leq k$. 
\end{Prop}

\noindent {\itshape Proof.} (i), (ii), and (iii) are clear by Proposition 
\ref{commformbeta}. We shall show (iv). 
Suppose that $k \neq \nu$. We proceed 
by induction on $k-j$. 
Since $e_{\beta_{k+1}}^{(c)}e_{\beta_k}^{(a_k)} -e_{\beta_k}^{(a_k)} e_{\beta_{k+1}}^{(c)}=0$
by the remark of Proposition \ref{commformbeta}, it is clear for $k-j=0$. So 
we assume that $k-j>0$.  We have 
\begin{align*}
e_{\beta_{k+1}}^{(c)}e_{\beta_j}^{(a_j)} \cdots e_{\beta_k}^{(a_k)} - 
e_{\beta_j}^{(a_j)} \cdots e_{\beta_k}^{(a_k)} e_{\beta_{k+1}}^{(c)} 
&= \left( e_{\beta_{k+1}}^{(c)}e_{\beta_j}^{(a_j)} - e_{\beta_j}^{(a_j)} e_{\beta_{k+1}}^{(c)}\right)
e_{\beta_{j+1}}^{(a_{j+1})} \cdots  e_{\beta_k}^{(a_k)} \\
& \ \  +  e_{\beta_j}^{(a_j)} \left(e_{\beta_{k+1}}^{(c)}e_{\beta_{j+1}}^{(a_{j+1})} \cdots 
e_{\beta_k}^{(a_k)} - e_{\beta_{j+1}}^{(a_{j+1})} \cdots e_{\beta_k}^{(a_k)} 
e_{\beta_{k+1}}^{(c)}\right)
\end{align*}
in $\mathcal{U}$. Now Proposition \ref{commformbeta} implies that  
\[e_{\beta_{k+1}}^{(c)}e_{\beta_j}^{(a_j)} - e_{\beta_j}^{(a_j)} e_{\beta_{k+1}}^{(c)}
= \sum_{{\bm b}=(b_j, \dots, b_{k+1})} \xi_1({\bm b}) 
e_{\beta_{j}}^{(b_{j})} \cdots  e_{\beta_{k+1}}^{(b_{k+1})}\] 
in $\mathcal{U}$, where $ \xi_1({\bm b}) \in \mathbb{F}_p$ and  $b_{k+1}< c$, $b_j < a_j$, 
$\sum_{i=j}^k b_i \leq a_j$, and $\sum_{i=j+1}^{k+1} b_i \leq c$ for each ${\bm b}$ 
with $\xi_1({\bm b}) \neq 0$ (in particular, note that $b_i \in \mathcal{N}_r$ 
for $j \leq i \leq k$). In turn, by induction we have 
\[e_{\beta_{k+1}}^{(b_{k+1})}e_{\beta_{j+1}}^{(a_{j+1})} \cdots e_{\beta_k}^{(a_k)} 
= e_{\beta_{j+1}}^{(a_{j+1})} \cdots e_{\beta_k}^{(a_k)} e_{\beta_{k+1}}^{(b_{k+1})} +
\sum_{{\bm d}=(d_{j+1}, \dots, d_{k+1})} \xi_2({\bm d})
e_{\beta_{j+1}}^{(d_{j+1})} \cdots e_{\beta_{k+1}}^{(d_{k+1})}\]
in $\mathcal{U}$, where  $ \xi_2({\bm d}) \in \mathbb{F}_p$ and  $d_{k+1}< b_{k+1}$ and 
$d_i \in \mathcal{N}_r$ for $j+1 \leq i \leq k$. Now we have
\begin{align*}
\lefteqn{\left( e_{\beta_{k+1}}^{(c)}e_{\beta_j}^{(a_j)} - e_{\beta_j}^{(a_j)} 
e_{\beta_{k+1}}^{(c)}\right)e_{\beta_{j+1}}^{(a_{j+1})} \cdots  e_{\beta_k}^{(a_k)} } \\
&= \sum_{{\bm b}=(b_j, \dots, b_{k+1})} \xi_1({\bm b}) 
e_{\beta_{j}}^{(b_{j})} \cdots  e_{\beta_{k+1}}^{(b_{k+1})} 
e_{\beta_{j+1}}^{(a_{j+1})} \cdots  e_{\beta_{k}}^{(a_{k})} \\
&=  \sum_{{\bm b}} \xi_1({\bm b}) 
e_{\beta_{j}}^{(b_{j})} \cdots  e_{\beta_{k}}^{(b_{k})} 
\left( e_{\beta_{j+1}}^{(a_{j+1})} \cdots  e_{\beta_{k}}^{(a_{k})}  e_{\beta_{k+1}}^{(b_{k+1})} 
+ \sum_{{\bm d}=(d_{j+1}, \dots, d_{k+1})} \xi_2({\bm d})
e_{\beta_{j+1}}^{(d_{j+1})} \cdots e_{\beta_{k+1}}^{(d_{k+1})}\right)
\end{align*} 
in $\mathcal{U}$. By (iii), the element is an 
$\mathbb{F}_p$-linear combination of elements of the form 
$e_{\beta_{j}}^{(c_{j})} \cdots  e_{\beta_{k+1}}^{(c_{k+1})}$ satisfying 
$c_{k+1}< c$ and $c_i \in \mathcal{N}_r$ for $j \leq i \leq k$. 

On the other hand, consider the element 
\[e_{\beta_j}^{(a_j)} \left(e_{\beta_{k+1}}^{(c)}e_{\beta_{j+1}}^{(a_{j+1})} \cdots 
e_{\beta_k}^{(a_k)} - e_{\beta_{j+1}}^{(a_{j+1})} \cdots e_{\beta_k}^{(a_k)} 
e_{\beta_{k+1}}^{(c)}\right).\]
By induction, we have 
\[e_{\beta_{k+1}}^{(c)}e_{\beta_{j+1}}^{(a_{j+1})} \cdots 
e_{\beta_k}^{(a_k)} - e_{\beta_{j+1}}^{(a_{j+1})} \cdots e_{\beta_k}^{(a_k)} 
e_{\beta_{k+1}}^{(c)} 
= \sum_{{\bm b}=(b_{j+1}, \dots, b_{k+1})} \xi_3({\bm b}) 
e_{\beta_{j+1}}^{(b_{j+1})} \cdots  e_{\beta_{k+1}}^{(b_{k+1})} \]
in $\mathcal{U}$, where  $ \xi_3({\bm b}) \in \mathbb{F}_p$ and  $b_{k+1}< c$ and 
$b_i \in \mathcal{N}_r$ for $j+1 \leq i \leq k$. 
Therefore, the result for $k \neq \nu$ follows. The proof for $j \neq 1$ is similar. 
$\square$ \\ \\

\section{Some linear isomorphisms for $\mathcal{U}^+$}

The aim in this section is to give several $\mathbb{F}_p$-linear isomorphisms induced by 
multiplication in $\mathcal{U}^+$. Throughout this section, let $r$ be a fixed positive integer. 

We need some propositions and lemmas to  prove  main results in this section. \\ 

\begin{Prop}\label{fr}
For $\beta \in \Phi^+$, the element ${\rm Fr}'^r(e_{\beta})-e_{\beta}^{(p^r)}$ 
in $\mathcal{U}$ lies in $\mathcal{U}_{r}^+$. 
\end{Prop}

\noindent {\itshape Proof.} 
We proceed by induction on ${\rm ht}(\beta)$. It is clear for ${\rm ht}(\beta)=1$ (i.e. $\beta$ is simple), so we assume that 
${\rm ht}(\beta) \geq 2$. 

We choose $\alpha \in \Delta$ such that $\beta-\alpha \in \Phi^+$ and let $k$ be a unique positive integer 
such that $\beta- k\alpha \in \Phi^+$ and  $\beta- (k+1)\alpha \not\in \Phi^+$. Then note that  
$[e_{\alpha}, e_{\beta-\alpha}] =  k c e_{\beta}$ 
in $\mathfrak{g}_{\mathbb{Z}}$, where $c=c_{\alpha, \beta-\alpha} \in \{ \pm 1\}$. 
Note also that $k \in \{1,2,3\}$. 

Suppose that $k \neq p$. By induction, the element $z_1 = {\rm Fr}'^r(e_{\beta-\alpha})-e_{\beta-\alpha}^{(p^r)}$ in 
$\mathcal{U}$ lies in $\mathcal{U}_{r}^+$. Then  we have 
\begin{align*}
\lefteqn{{\rm Fr}'^r(e_{\beta})-e_{\beta}^{(p^r)} 
= {\rm Fr}'^r \left( k^{-1} c (e_{\alpha}e_{\beta-\alpha}-e_{\beta-\alpha}e_{\alpha})\right)
-e_{\beta}^{(p^r)}} \\
&=k^{-1}  c 
\left( e_{\alpha}^{(p^r)} {\rm Fr}'^r(e_{\beta-\alpha}) - 
{\rm Fr}'^r(e_{\beta-\alpha})e_{\alpha}^{(p^r)} \right) -e_{\beta}^{(p^r)} \\
&=  k^{-1} c
\left( e_{\alpha}^{(p^r)} e_{\beta-\alpha}^{(p^r)} - 
e_{\beta-\alpha}^{(p^r)}e_{\alpha}^{(p^r)} -kc e_{\beta}^{(p^r)} \right) 
+ k^{-1} c (e_{\alpha}^{(p^r)} z_1 - z_1 e_{\alpha}^{(p^r)})
\end{align*}
in $\mathcal{U}$, where $k^{-1}$ denotes the inverse  of  $k$ in $\mathbb{F}_p$. 
Now taking $n=1$ and $\gamma = \beta-\alpha$ in Proposition \ref{commformpr1}, 
we see that   the element 
$e_{\alpha}^{(p^r)} e_{\beta-\alpha}^{(p^r)} - 
e_{\beta-\alpha}^{(p^r)}e_{\alpha}^{(p^r)} -kc e_{\beta}^{(p^r)}$ lies in $\mathcal{U}_r$ 
(hence in $\mathcal{U}_r^+$). On the other hand, 
taking $n=1$ and $z=z_1$ in Proposition \ref{commformpr2}, we see that the element  
$e_{\alpha}^{(p^r)} z_1 - z_1 e_{\alpha}^{(p^r)}$ lies in $\mathcal{U}_{r}$ 
(hence in $\mathcal{U}_r^+$). 
Therefore, the element 
${\rm Fr}'^r(e_{\beta})-e_{\beta}^{(p^r)}$  lies in $\mathcal{U}_{r}^+$. 

Suppose that $k=p=2$. Then
$[e_{\alpha}, e_{\beta-2\alpha}]= c_1 e_{\beta-\alpha}$
and 
$[e_{\alpha}, e_{\beta-\alpha}]= 2c_2 e_{\beta}$ 
in $\mathfrak{g}_{\mathbb{Z}}$, where 
$c_1=c_{\alpha, \beta-2\alpha} \in \{ \pm 1 \}$ and 
$c_2=c_{\alpha, \beta-\alpha} \in \{ \pm 1 \}$. Then we have 
\[e_{\beta}=c_1c_2 (e_{\alpha}^{(2)} e_{\beta-2\alpha} 
- e_{\alpha} e_{\beta-2\alpha} e_{\alpha} +e_{\beta-2\alpha} e_{\alpha}^{(2)})\]
in $\mathcal{U}_{\mathbb{Z}}$. 
By induction, the element 
$z_2 = {\rm Fr}'^r(e_{\beta-2\alpha})-e_{\beta-2\alpha}^{(2^r)}$ in 
$\mathcal{U}$ lies in $\mathcal{U}_{r}^+$. Then  we have 
\begin{align*}
\lefteqn{{\rm Fr}'^r(e_{\beta})-e_{\beta}^{(2^r)}} \\
&= c_1 c_2 \left( e_{\alpha}^{(2^{r+1})} {\rm Fr}'^r(e_{\beta-2\alpha}) - 
e_{\alpha}^{(2^{r})}{\rm Fr}'^r(e_{\beta-2\alpha})e_{\alpha}^{(2^{r})}+
{\rm Fr}'^r(e_{\beta-2\alpha})e_{\alpha}^{(2^{r+1})}\right) - e_{\beta}^{(2^r)} \\
&= \left( e_{\alpha}^{(2^{r+1})} e_{\beta-2\alpha}^{(2^r)} - 
e_{\alpha}^{(2^{r})} e_{\beta-2\alpha}^{(2^r)} e_{\alpha}^{(2^{r})}+
e_{\beta-2\alpha}^{(2^r)} e_{\alpha}^{(2^{r+1})} - e_{\beta}^{(2^r)}\right) \\
& \ \ \ +\left( e_{\alpha}^{(2^{r+1})} z_2- 
e_{\alpha}^{(2^{r})} z_2 e_{\alpha}^{(2^{r})}+
z_2 e_{\alpha}^{(2^{r+1})}\right)
\end{align*}
in $\mathcal{U}$ (since $p=2$, note that $c_1=c_2=1$ in $\mathbb{F}_p$). 
Now taking $p=2$, $n=2$, and $\gamma=\beta- 2\alpha$ in 
Proposition \ref{commformpr1}, we see that 
$\alpha$ and $\gamma$ are of the case (D) or (F) and that 
the element 
\[
e_{\alpha}^{(2^{r+1})} e_{\beta-2\alpha}^{(2^r)} - 
e_{\alpha}^{(2^{r})} e_{\beta-2\alpha}^{(2^r)} e_{\alpha}^{(2^{r})}+
e_{\beta-2\alpha}^{(2^r)} e_{\alpha}^{(2^{r+1})} - e_{\beta}^{(2^r)}
\]
lies in 
$\mathcal{U}_r$ (hence in $\mathcal{U}_r^+$). On the other hand, taking 
$p=2$, $n=2$, and $z=z_2$ in 
Proposition \ref{commformpr2}, we see that 
the element 
$e_{\alpha}^{(2^{r+1})} z_2- 
e_{\alpha}^{(2^{r})} z_2 e_{\alpha}^{(2^{r})}+
z_2 e_{\alpha}^{(2^{r+1})}$  lies in 
$\mathcal{U}_r$ (hence in $\mathcal{U}_r^+$).  
Therefore, the element  ${\rm Fr}'^r(e_{\beta})-e_{\beta}^{(2^r)}$ lies in $\mathcal{U}_{r}^+$. 

Finally, suppose that $k=p=3$. This case can occur only when 
$\Phi$ is of type ${\rm G}_2$, $\beta = 3\alpha_1+\alpha_2$, and $\alpha=\alpha_1$. 
Note that 
\[e_{1112}=e_1^{(3)} e_2 -e_1^{(2)} e_2 e_1 +e_1 e_2e_1^{(2)}-e_2e_1^{(3)}\]
in $\mathcal{U}_{\mathbb{Z}}$ (see Section 2 for the notation). 
Now taking $p=3$, $n=3$, $\alpha=\alpha_1$, and $\gamma=\beta-3\alpha=\alpha_2$ in 
Proposition \ref{commformpr1}, we see  
that $\alpha$ and $\gamma$ are of the case (F) and that the element  
\begin{align*}
\lefteqn{{\rm Fr}'^r(e_{1112})- e_{1112}^{(3^r)}} \\
&= e_1^{(3^{r+1})} e_2^{(3^r)} -e_1^{(2\cdot 3^r)} e_2^{(3^r)} e_1^{(3^r)} 
+e_1^{(3^r)} e_2^{(3^r)}e_1^{(2\cdot 3^r)}-e_2^{(3^r)}e_1^{(3^{r+1})} -e_{1112}^{(3^r)}
\end{align*} 
in $\mathcal{U}$ lies in $\mathcal{U}_{r}^+$. $\square$ \\
\ 

For $c \in \mathbb{Z}_{\geq 0}$, if we write $c=c' + p^r c''$ uniquely 
with $c' \in \mathcal{N}_r$ and 
$c'' \in \mathbb{Z}_{\geq 0}$, then we denote $c'$ by $r_{p,r}(c)$ and 
$c''$ by $q_{p,r}(c)$. 
Clearly we have $q_{p,r}(c)=0$ if and only if $c \in \mathcal{N}_r$. 
Moreover, for $n \in \mathbb{Z}_{> 0}$ and 
${\bm c}=(c_1, \dots, c_n) \in (\mathbb{Z}_{\geq 0})^n$, set 
\[q_{p,r}({\bm c})=(q_{p,r}(c_1), \dots, q_{p,r}(c_n))\]
and 
\[r_{p,r}({\bm c})=(r_{p,r}(c_1), \dots, r_{p,r}(c_n)).\] 

For ${\bm a}=(a_{j}, \dots, a_{k}) \in (\mathbb{Z}_{\geq 0})^{k-j+1}$ with 
$1 \leq j \leq k \leq \nu$, we shall often 
write the element $e_{\beta_{j}}^{(a_j)} \cdots e_{\beta_{k}}^{(a_k)}$ 
(resp. $e_{-\beta_{j}}^{(a_j)} \cdots e_{-\beta_{k}}^{(a_k)}$) of 
$\mathcal{U}_{\mathbb{Z}}$ or $\mathcal{U}$ as ${\bm e}^{({\bm a})}$ 
(resp. ${\bm f}^{({\bm a})}$) for simplicity. \\ 
\ 

\begin{Lem}\label{fact1} 
Let ${\bm a}=(a_1, \dots, a_{\nu}) \in (\mathbb{Z}_{\geq 0})^{\nu}$. Then the following hold. \\ 

\noindent {\rm (i)} Suppose that there exists an integer $k$ with $1 \leq k \leq \nu$ 
such that $q_{p,r}(a_k)<p-1$ and $q_{p,r}(a_i)=0$ whenever $k+1 \leq i \leq \nu$. Then 
we have 
\[{\bm e}^{({\bm a})} e_{\beta_{k}}^{(p^r)} 
= (q_{p,r}(a_k)+1) {\bm e}^{(\widetilde{\bm a})} + 
\sum_{{\bm c}=(c_1, \dots, c_{\nu})} \xi({\bm c}) {\bm e}^{({\bm c})}\]
in $\mathcal{U}$, 
where $\widetilde{\bm a}=(a_1, \dots, a_{k-1}, a_k+p^r, a_{k+1}, \dots, a_{\nu})$, 
$\xi({\bm c}) \in \mathbb{F}_p$, and each ${\bm c}$ with 
$\xi({\bm c}) \neq 0$ satisfies $q_{p,r}({\bm c})= q_{p,r}({\bm a})$ {\rm (}i.e. 
$q_{p,r}(c_i)= q_{p,r}(a_i)$ for $1 \leq i \leq \nu${\rm )}. \\

\noindent {\rm (ii)} For $z \in \mathcal{U}_r^+$, if we write the element 
${\bm e}^{({\bm a})} z$ in $\mathcal{U}^+$ as the form 
$\sum_{{\bm c}=(c_1, \dots, c_{\nu})} \xi({\bm c}) {\bm e}^{({\bm c})}$ with 
$\xi({\bm c}) \in \mathbb{F}_p$, then each ${\bm c}$ with $\xi({\bm c}) \neq 0$ satisfies 
$q_{p,r}(c_i) \leq q_{p,r}(a_i)$ for $1 \leq i \leq \nu$. 
\end{Lem}

\noindent {\itshape Proof.} (i) Since $q_{p,r}(a_k)< p-1$, by Proposition  
\ref{basicformulas} (i) and Proposition \ref{binomial1} we have 
\[e_{\beta_k}^{(a_k)} e_{\beta_k}^{(p^r)} = 
{a_k +p^r \choose p^r} e_{\beta_k}^{(a_k+p^r)} = (q_{p,r}(a_k)+1)  e_{\beta_k}^{(a_k+p^r)} \] 
in $\mathcal{U}$. So it is clear when $k=\nu$.  
Now assume that $k < \nu$. 
Applying Proposition \ref{commformbeta} to 
the element $e_{\beta_{k+1}}^{(a_{k+1})} \cdots e_{\beta_{\nu}}^{(a_{\nu})} 
e_{\beta_{k}}^{(p^r)}$ repeatedly, we have 
\[e_{\beta_{k+1}}^{(a_{k+1})} \cdots e_{\beta_{\nu}}^{(a_{\nu})} 
e_{\beta_{k}}^{(p^r)}
= e_{\beta_{k}}^{(p^r)} e_{\beta_{k+1}}^{(a_{k+1})} \cdots e_{\beta_{\nu}}^{(a_{\nu})} 
+ \sum_{{\bm c}=(c_k, \dots, c_{\nu})} \xi'({\bm c}) {\bm e}^{({\bm c})}\]
in $\mathcal{U}$, where $\xi'({\bm c}) \in \mathbb{F}_p$ and each ${\bm c}$ with 
$\xi'({\bm c}) \neq 0$ satisfies 
$c_i \in \mathcal{N}_r$ for $k \leq i \leq \nu$. Thus we have 
\begin{align*}
\lefteqn{{\bm e}^{({\bm a})} e_{\beta_k}^{(p^r)} =
e_{\beta_1}^{(a_1)} \cdots e_{\beta_k}^{(a_k)} \cdot e_{\beta_{k+1}}^{(a_{k+1})} \cdots 
e_{\beta_{\nu}}^{(a_{\nu})} e_{\beta_k}^{(p^r)} } \\
&= e_{\beta_1}^{(a_1)} \cdots e_{\beta_k}^{(a_k)} 
\left( e_{\beta_k}^{(p^r)}e_{\beta_{k+1}}^{(a_{k+1})} \cdots 
e_{\beta_{\nu}}^{(a_{\nu})} + \sum_{{\bm c}=(c_k, \dots, c_{\nu})} \xi'({\bm c})  
{\bm e}^{({\bm c})}\right) \\
&= (q_{p,r}(a_k)+1) {\bm e}^{(\widetilde{\bm a})} 
+\sum_{{\bm c}=(c_k, \dots, c_{\nu})} \xi'({\bm c}) e_{\beta_1}^{(a_1)}  \cdots 
e_{\beta_{k-1}}^{(a_{k-1})} e_{\beta_k}^{(a_k)}  e_{\beta_k}^{(c_k)} \cdots 
e_{\beta_{\nu}}^{(c_{\nu})}   \\
&= (q_{p,r}(a_k)+1) {\bm e}^{(\widetilde{\bm a})} 
+\sum_{{\bm c}=(c_k, \dots, c_{\nu})} \xi'({\bm c}) {a_k+c_k \choose c_k} 
e_{\beta_1}^{(a_1)}  \cdots 
e_{\beta_{k-1}}^{(a_{k-1})} e_{\beta_k}^{(a_k+c_k)}  e_{\beta_{k+1}}^{(c_{k+1})} \cdots 
e_{\beta_{\nu}}^{(c_{\nu})}  
\end{align*}
in $\mathcal{U}$. 
Now if $q_{p,r}(a_k)< q_{p,r}(a_k+c_k)$, then $r_{p,r}(a_k+c_k)< c_k$ and hence 
${a_k+c_k \choose c_k} =0$ in $\mathbb{F}_p$. 
So in the last sum of the above equality a term with nonzero coefficient must satisfy 
$q_{p,r}(a_k) =q_{p,r}(a_k+c_k)$. Moreover, we have 
$q_{p,r}(a_i)=0=q_{p,r}(c_i)$ for $k+1 \leq i \leq \nu$. Therefore, (i) follows. 

(ii) More generally, we shall prove the following lemma: \\

\begin{Lem}\label{fact1'}
Suppose that $j,k \in \{1,2, \dots, \nu \}$, $j-1 \leq k$, and $z \in \mathcal{U}_r^+$. 
Then if we write  the element $\left( \prod_{i=j}^{k} e_{\beta_i}^{(a_i)} \right) z$ in 
$\mathcal{U}^+$ with 
$a_i \in \mathbb{Z}_{\geq 0}$ as the form 
$\sum_{{\bm c}=(c_1, \dots, c_{\nu})} \xi({\bm c}) {\bm e}^{({\bm c})}$ with 
$\xi({\bm c}) \in \mathbb{F}_p$, then each ${\bm c}$ with $\xi({\bm c}) \neq 0$ satisfies 
$q_{p,r}(c_i) \leq q_{p,r}(a_i)$ for $j \leq i \leq k$ and $q_{p,r}(c_i) =0$ for 
$1 \leq i <j$ or $k < i \leq \nu$.    
\end{Lem}

\noindent {\itshape Proof.} We may assume that 
$z=\prod_{i=1}^{\nu} e_{\beta_i}^{(b_i)}$ with $b_i \in \mathcal{N}_r$ for each $i$. 
We proceed by induction on $k-j+1$. 
If $k-j+1=0$, then $\prod_{i=j}^{k} e_{\beta_i}^{(a_i)}=1$ and the result is clear. 
So from now on assume that $k \geq j$. By induction we have 
\[ \left( \prod_{i=j+1}^{k} e_{\beta_i}^{(a_i)} \right) z 
= \sum_{{\bm d}=(d_1, \dots, d_{\nu})} \xi_1({\bm d}) {\bm e}^{({\bm d})}\]
in $\mathcal{U}$, where $\xi_1({\bm d}) \in \mathbb{F}_p$ and each ${\bm d}$ with 
$\xi_1({\bm d})\neq 0$ satisfies 
$q_{p,r}(d_i) \leq q_{p,r}(a_i)$ for $j+1 \leq i \leq k$ and $q_{p,r}(d_i)=0$ for 
$1 \leq i < j+1$ or $k < i \leq \nu$. 
Thus, to prove the lemma, it suffices to show that for each ${\bm d}$ with 
$\xi_1({\bm d})\neq 0$, if we write the element $e_{\beta_j}^{(a_j)} {\bm e}^{({\bm d})}$ 
in $\mathcal{U}$ as 
the form $\sum_{{\bm c}'=(c_1', \dots, c_{\nu}')} \xi'({\bm c}') {\bm e}^{({\bm c}')}$ 
with $\xi'({\bm c}') \in \mathbb{F}_p$, then each ${\bm c}'$ with $\xi'({\bm c}') \neq 0$ 
satisfies $q_{p,r}(c_i') \leq q_{p,r}(a_i)$ for $j \leq i \leq k$ and $q_{p,r}(c_i')=0$ for 
$1 \leq i < j$ or $k < i \leq \nu$. 
It is clear when $a_j=0$, so we may assume that $a_j>0$. 
Let ${\bm d}=(d_1, \dots, d_{\nu}) 
\in (\mathbb{Z}_{\geq 0})^{\nu}$ be any $\nu$-tuple with $\xi_1({\bm d}) \neq 0$. 
Then recall that ${\bm d}$ satisfies $q_{p,r}(d_i) \leq q_{p,r}(a_i)$ for $j+1 \leq i \leq k$ and 
$d_i \in \mathcal{N}_r$ for $1 \leq i < j+1$ or $k<i \leq \nu$. 
By Proposition \ref{subring} (iv), we have 
\[e_{\beta_j}^{(a_j)} \prod_{i=1}^{j-1}e_{\beta_i}^{(d_i)} = 
\left( \prod_{i=1}^{j-1}e_{\beta_i}^{(d_i)} \right) e_{\beta_j}^{(a_j)} 
+ \sum_{{\bm d}'=(d_1', \dots, d_{j}')} \xi_2({\bm d}') {\bm e}^{({\bm d}')}\]
in $\mathcal{U}$, where $\xi_2({\bm d}') \in \mathbb{F}_p$ and each ${\bm d}'$ with 
$\xi_2({\bm d}') \neq 0$ satisfies $d_j' < a_j$ 
(hence $q_{p,r}(d_j') \leq q_{p,r}(a_j)$) and $d_i' \in \mathcal{N}_r$  
(i.e. $q_{p,r}(d_i')=0$) for 
$1 \leq i \leq j-1$. Then we have 
\begin{align*}
e_{\beta_j}^{(a_j)} {\bm e}^{({\bm d})}
&=\left( \left( \prod_{i=1}^{j-1}e_{\beta_i}^{(d_i)} \right) e_{\beta_j}^{(a_j)} 
+ \sum_{{\bm d}'=(d_1', \dots, d_{j}')} \xi_2({\bm d}') {\bm e}^{({\bm d}')}\right) 
\prod_{i=j}^{\nu} e_{\beta_i}^{(d_i)}   \\
&= {a_j +d_j \choose d_j}  \left( \prod_{i=1}^{j-1}e_{\beta_i}^{(d_i)} \right) 
e_{\beta_j}^{(a_j+d_j)} \prod_{i=j+1}^{\nu} e_{\beta_i}^{(d_i)} \\
& \ \ \ + \sum_{{\bm d}'=(d_1', \dots, d_{j}')} \xi_2({\bm d}')  
{d_j' +d_j \choose d_j}  \left( \prod_{i=1}^{j-1}e_{\beta_i}^{(d_i')} \right) e_{\beta_j}^{(d_j'+d_j)} \prod_{i=j+1}^{\nu} e_{\beta_i}^{(d_i)}
\end{align*}
in $\mathcal{U}$. Since $q_{p,r}(d_j)=0$,  
if ${a_j +d_j \choose d_j} \neq 0$ in $\mathbb{F}_p$, then 
$q_{p,r}(a_j+d_j)=q_{p,r}(a_j)$ and similarly,  if  ${d_j' +d_j \choose d_j} \neq 0$  
in $\mathbb{F}_p$, then $q_{p,r}(d_j'+d_j)=q_{p,r}(d_j')\ (\leq q_{p,r}(a_j))$. 
Finally, recall that $q_{p,r}(d_i')=0$ for $1 \leq i < j$, 
$q_{p,r}(d_i) \leq q_{p,r}(a_i)$ for $j+1 \leq i \leq k$, and $q_{p,r}(d_i)=0$ for 
$1 \leq i < j$ or $k < i \leq \nu$.  Therefore, 
the lemma follows. $\square$ \\

Now Lemma \ref{fact1} (ii) follows by taking $j=1$ and $k=\nu$ in Lemma \ref{fact1'}. 
Therefore, we complete the proof of Lemma \ref{fact1}. $\square$ \\

\begin{Prop}\label{fact2}
Suppose that  ${\bm a}=(a_1, \dots, a_{\nu}) \in (\mathcal{N}_r)^{\nu}$ 
 and that ${\bm b}=(b_1, \dots, b_{k}) \in (\mathcal{N}_1)^{k}$ 
with $1 \leq k \leq \nu$. Then  we have 
\[ {\bm e}^{(\bm a)} {\rm Fr}'^r \left( {\bm e}^{({\bm b})}\right) =
\left( \prod_{i=1}^k e_{\beta_i}^{(a_i+p^r b_i)}\right) \prod_{i=k+1}^{\nu} e_{\beta_i}^{(a_i)}
+ \sum_{{\bm c}=(c_1, \dots, c_{\nu})} \xi({\bm c}) {\bm e}^{({\bm c})}\]
in $\mathcal{U}$, where $\xi({\bm c}) \in \mathbb{F}_p$ and each 
${\bm c}$ with $\xi({\bm c}) \neq 0$ satisfies 
\[(q_{p,r}(c_1), \dots, q_{p,r}(c_k)) \neq 
(b_1, \dots, b_k)\]
in $(\mathbb{Z}_{\geq 0})^k$, $q_{p,r}(c_i) \leq b_i$ for $1 \leq i \leq k$, 
and $q_{p,r}(c_i)=0$ for $k+1 \leq i \leq \nu$. 
\end{Prop}

\noindent {\itshape Proof.} 
We proceed by induction on $n=\sum_{i=1}^k b_i$. 
It is clear for $n=0$, so assume that $n>0$. We may also assume that $b_k \neq 0$. 
Set 
\[
\widetilde{\bm b}=(b_1, \dots, b_{k-1}, b_k-1)
\]
and 
\[
\widehat{\bm a}=\left( a_1+p^r b_1, \dots, a_{k-1}+p^r b_{k-1},a_k+p^r(b_k-1), 
a_{k+1}, \dots, a_{\nu}\right)
.\] 
We have  
\[{\bm e}^{({\bm a})} {\rm Fr}'^r \left( {\bm e}^{({\bm b})}\right) = 
b_k^{-1} {\bm e}^{({\bm a})} 
{\rm Fr}'^r \left( {\bm e}^{(\widetilde{\bm b})} e_{\beta_k}\right) =
b_k^{-1} {\bm e}^{({\bm a})} 
{\rm Fr}'^r \left( {\bm e}^{(\widetilde{\bm b})}\right)  
{\rm Fr}'^r \left(e_{\beta_k}\right) \]
in $\mathcal{U}$, where $b_k^{-1}$ is the inverse of $b_k$ in $\mathbb{F}_p$. 
By induction, we have 
\[{\bm e}^{({\bm a})} {\rm Fr}'^r \left( {\bm e}^{(\widetilde{\bm b})}\right) = 
{\bm e}^{(\widehat{\bm a})} 
+ \sum_{{\bm c}=(c_1, \dots, c_{\nu})} \xi_1({\bm c}) {\bm e}^{({\bm c})}\]
in $\mathcal{U}$, where $ \xi_1({\bm c}) \in \mathbb{F}_p$ and 
each ${\bm c}$ with $\xi_1({\bm c}) \neq 0$ 
satisfies 
\[
(q_{p,r}(c_1), \dots, q_{p,r}(c_{k-1}), q_{p,r}(c_{k})) \neq 
(b_1, \dots, b_{k-1}, b_k-1)
\]
in $(\mathbb{Z}_{\geq 0})^k$, $q_{p,r}(c_i) \leq b_i$ for $1 \leq i \leq k-1$, $q_{p,r}(c_k) \leq b_k-1$, 
 and $q_{p,r}(c_i)=0$ for $k+1 \leq i \leq \nu$. Moreover, by Proposition 
\ref{fr} there exists $z \in \mathcal{U}_r^+$ such that 
${\rm Fr}'^r(e_{\beta_k})= e_{\beta_k}^{(p^r)} +z$. Then we have 
\begin{align*}
\lefteqn{{\bm e}^{(\bm a)} {\rm Fr}'^r \left( {\bm e}^{({\bm b})}\right) =
b_k^{-1} \left( {\bm e}^{(\widehat{\bm a})} 
+ \sum_{{\bm c}=(c_1, \dots, c_{\nu})} \xi_1({\bm c}) {\bm e}^{({\bm c})} \right) 
\left( e_{\beta_k}^{(p^r)} +z \right)} \\
&= b_k^{-1} \left( {\bm e}^{(\widehat{\bm a})} e_{\beta_k}^{(p^r)}+ 
\sum_{\bm c} \xi_1({\bm c}) {\bm e}^{({\bm c})} e_{\beta_k}^{(p^r)}+
{\bm e}^{(\widehat{\bm a})} z + \sum_{\bm c} \xi_1({\bm c}) {\bm e}^{({\bm c})} z \right)
\end{align*}
in $\mathcal{U}$. Applying Lemma \ref{fact1} (i) to the element 
${\bm e}^{(\widehat{\bm a})} e_{\beta_k}^{(p^r)}$ we have 
\[
{\bm e}^{(\widehat{\bm a})} e_{\beta_k}^{(p^r)} = 
b_k {\bm e}^{(\widehat{\bm a}')} 
+\sum_{{\bm d}=(d_1, \dots, d_{\nu})} \xi_2({\bm d}) {\bm e}^{({\bm d})}
\]
in $\mathcal{U}$, where $\xi_2({\bm d}) \in \mathbb{F}_p$,   
\[\widehat{\bm a}'=\left( a_1+p^r b_1, \dots, a_{k-1}+p^r b_{k-1},a_k+p^rb_k, 
a_{k+1}, \dots, a_{\nu}\right),\]
and each ${\bm d}$ with $\xi_2({\bm d}) \neq 0$ satisfies 
$q_{p,r}({\bm d})=q_{p,r}(\widehat{\bm a})$. Then clearly we have 
\[
(q_{p,r}(d_1), \dots, q_{p,r}(d_{k-1}), q_{p,r}(d_{k})) = (b_1, \dots, b_{k-1}, b_k-1) \neq 
(b_1, \dots, b_{k-1}, b_k)
\] 
in $(\mathbb{Z}_{\geq 0})^k$, 
$q_{p,r}(d_i) \leq b_i$ for $1 \leq i \leq k$, and 
$q_{p,r}(d_i)=0$ for $k+1 \leq i \leq \nu$. On the other hand, applying 
Lemma \ref{fact1} (i) to the element $ {\bm e}^{({\bm c})} e_{\beta_k}^{(p^r)}$ and 
Lemma \ref{fact1} (ii) to the elements ${\bm e}^{(\widehat{\bm a})} z$ and 
${\bm e}^{({\bm c})} z$, we see that these elements are written as the form 
$\sum_{{\bm d}=(d_1, \dots, d_{\nu})} \xi_3({\bm d}) {\bm e}^{({\bm d})}$, 
where $\xi_3({\bm d}) \in \mathbb{F}_p$ and each ${\bm d}$ with $\xi_3({\bm d}) \neq 0$ 
satisfies 
\[
(q_{p,r}(d_1), \dots, q_{p,r}(d_{k})) \neq 
(b_1, \dots, b_k)
\] 
in $(\mathbb{Z}_{\geq 0})^k$, 
$q_{p,r}(d_i) \leq b_i$ for $1 \leq i \leq k$, and 
$q_{p,r}(d_i)=0$ for $k+1 \leq i \leq \nu$. 
Now the proposition follows. 
$\square$ \\

Now we are ready to describe main results in this section. 
\ \\  

\begin{The}\label{mainthm1}
Let $n \in \mathbb{Z}_{> 0}$. Then the multiplication on $\mathcal{U}^+$ induces 
the following two $\mathbb{F}_p$-linear isomorphisms:  
\[
\mathcal{U}_r^+ \otimes_{\mathbb{F}_p} {\rm Fr}'^r\left( \mathcal{U}_n^+\right) 
\rightarrow \mathcal{U}_{r+n}^+,\ \ \ 
\bigotimes_{i=0}^{r-1} {\rm Fr}'^i \left( \mathcal{U}_1^+\right) 
\rightarrow \mathcal{U}_{r}^+.
\] 
\end{The}
\ 

Before the proof, we shall introduce some orderings in $(\mathbb{Z}_{\geq 0})^n$ 
for a fixed $n \in \mathbb{Z}_{> 0}$. 
First, we define the lexicographical order $\geq$ in 
$(\mathbb{Z}_{\geq 0})^n$ such that 
\[
(0,0,\dots, 0,1) < (0,0,\dots, 1,0) < \cdots < (0,1,\dots, 0,0)  < (1,0,\dots, 0,0).  
\]
Then we define another ordering  $\underset{\text{$(p,r)$}}{\geq}$ in 
$(\mathbb{Z}_{\geq 0})^n$ as follows: For ${\bm a}=(a_1, \dots, a_n), 
{\bm b}=(b_1, \dots, b_n) \in (\mathbb{Z}_{\geq 0})^n$, we write 
${\bm a} \underset{\text{$(p,r)$}}{\geq} {\bm b}$ if 
$q_{p,r}({\bm a}) > q_{p,r}({\bm b})$ or if $q_{p,r}({\bm a}) = q_{p,r}({\bm b})$ and 
$r_{p,r}({\bm a}) \geq r_{p,r}({\bm b})$. Note that both of the orderings in 
$(\mathbb{Z}_{\geq 0})^n$ are total. 
 \\ 

\noindent {\itshape Proof of Theorem \ref{mainthm1}.}  
To begin with, we shall show the result for the first map in the theorem. 
We proceed by induction on $n$. Suppose that $n=1$.  The elements 
${\bm e}^{({\bm a})} \otimes {\rm Fr}'^r \left( {\bm e}^{({\bm b})}\right)$ with 
${\bm a} =(a_1, \dots, a_{\nu}) \in (\mathcal{N}_r)^{\nu}$ and 
${\bm b} =(b_1, \dots, b_{\nu}) \in (\mathcal{N}_1)^{\nu}$ form an $\mathbb{F}_p$-basis 
of $\mathcal{U}_r^+ \otimes_{\mathbb{F}_p} {\rm Fr}'^r \left( \mathcal{U}_1^+\right)$ 
(see \cite[II, 1.12]{jantzenbook}). 
Since the multiplication map 
$\mathcal{U}_r^+ \otimes_{\mathbb{F}_p} {\rm Fr}'^r\left( \mathcal{U}_1^+\right) 
\rightarrow \mathcal{U}_{r+1}^+$ is $\mathbb{F}_p$-linear and 
\[{\rm dim}_{\mathbb{F}_p} \left( \mathcal{U}_r^+ \otimes_{\mathbb{F}_p}  
{\rm Fr}'^r \left( \mathcal{U}_1^+\right) \right)= p^{r \nu} \cdot p^{\nu} =p^{(r+1)\nu} 
={\rm dim}_{\mathbb{F}_p} \mathcal{U}_{r+1}^+,\]
it suffices to show that the map is injective, hence that the elements 
${\bm e}^{({\bm a})}  {\rm Fr}'^r \left( {\bm e}^{({\bm b})}\right)$ in 
$\mathcal{U}_{r+1}^+$ with 
${\bm a} =(a_1, \dots, a_{\nu}) \in (\mathcal{N}_r)^{\nu}$ and 
${\bm b} =(b_1, \dots, b_{\nu}) \in (\mathcal{N}_1)^{\nu}$ are linearly independent over 
$\mathbb{F}_p$. Consider a linear relation 
\[
\sum_{({\bm a}, {\bm b}) \in (\mathcal{N}_r)^{\nu} \times (\mathcal{N}_1)^{\nu}}
\eta({\bm a}, {\bm b}){\bm e}^{({\bm a})}  {\rm Fr}'^r \left( {\bm e}^{({\bm b})}\right) =0 
\eqno{(*)}
\]
with $\eta({\bm a}, {\bm b}) \in \mathbb{F}_p$. We need to show that 
$\eta({\bm a}, {\bm b})=0$ for all 
$({\bm a}, {\bm b}) \in (\mathcal{N}_r)^{\nu} \times (\mathcal{N}_1)^{\nu}$. 
Suppose that there exists a pair 
$({\bm a}, {\bm b}) \in (\mathcal{N}_r)^{\nu} \times (\mathcal{N}_1)^{\nu}$ such that 
$\eta({\bm a}, {\bm b}) \neq 0$. 
Then there is a unique pair $(\widetilde{\bm a}, \widetilde{\bm b} ) 
\in (\mathcal{N}_r)^{\nu} \times (\mathcal{N}_1)^{\nu}$ 
with 
$\widetilde{\bm a} = (\widetilde{a}_1, \dots, \widetilde{a}_{\nu})$ and 
$\widetilde{\bm b} = (\widetilde{b}_1, \dots, \widetilde{b}_{\nu})$
such that 
\[\widetilde{\bm a}+p^r \widetilde{\bm b} = 
(\widetilde{a}_1+p^r \widetilde{b}_1, \dots, 
\widetilde{a}_{\nu}+p^r \widetilde{b}_{\nu})\]
is the largest element with respect to the ordering 
$\underset{\text{$(p,r)$}}{\geq}$ 
among all ${\bm a}+p^r {\bm b}$ for various pairs $({\bm a}, {\bm b})$ 
satisfying $\eta({\bm a}, {\bm b}) \neq 0$. 
By Proposition \ref{fact2}, for each 
$({\bm a}, {\bm b}) \in (\mathcal{N}_r)^{\nu} \times (\mathcal{N}_1)^{\nu}$, we can write 
the element ${\bm e}^{({\bm a})} {\rm Fr}'^r \left( {\bm e}^{({\bm b})} \right)$ as 
\[{\bm e}^{({\bm a})} {\rm Fr}'^r \left( {\bm e}^{({\bm b})} \right)
={\bm e}^{({\bm a}+p^r{\bm b})}
+\sum_{{\bm c}=(c_1,\dots, c_{\nu})} \xi_1({\bm c}) {\bm e}^{({\bm c})},\]
where $\xi_1({\bm c}) \in \mathbb{F}_p$ and each ${\bm c}$ with $\xi_1({\bm c}) \neq 0$ 
satisfies ${\bm c} \underset{\text{$(p,r)$}}{<} {\bm a}+p^r {\bm b}$. 
Thus the left-hand side of the above linear relation $(*)$ can be written as 
\[\eta (\widetilde{\bm a}, \widetilde{\bm b}) 
{\bm e}^{(\widetilde{\bm a}+p^r\widetilde{\bm b})} + 
\sum_{{\bm c}=(c_1,\dots, c_{\nu})} \xi_2({\bm c}) {\bm e}^{({\bm c})},\]
where $\xi_2({\bm c}) \in \mathbb{F}_p$ and each ${\bm c}$ with $\xi_2({\bm c}) \neq 0$ 
satisfies ${\bm c} \underset{\text{$(p,r)$}}{<} \widetilde{\bm a}+p^r\widetilde{\bm b}$. 
Then $\mathbb{F}_p$-linear independence of the elements 
${\bm e}^{({\bm d})}$ with ${\bm d} \in (\mathbb{Z}_{\geq 0})^{\nu}$ implies that 
$\eta (\widetilde{\bm a}, \widetilde{\bm b}) =0$, which is contradiction. 
Therefore, we have shown the theorem for $n=1$. 

From now on, assume that $n \geq 2$. Using the result of the last paragraph, 
the fact that  ${\rm Fr}'$ is an 
$\mathbb{F}_p$-algebra homomorphism on $\mathcal{U}^+$, and induction on $n$, 
we obtain the following commutative diagram of well-defined $\mathbb{F}_p$-linear maps 
induced by multiplication: 
\[
  \begin{CD}
     {\mathcal{U}_r^+ \otimes_{\mathbb{F}_p} {\rm Fr}'^r
\left( \mathcal{U}_{n-1}^+ \right)
\otimes_{\mathbb{F}_p} {\rm Fr}'^{r+n-1}
\left( \mathcal{U}_1^+ \right)} @>>> 
{\mathcal{U}_r^+ \otimes_{\mathbb{F}_p} {\rm Fr}'^r
\left( \mathcal{U}_{n}^+ \right)} \\
  @VVV    @VVV \\
{\mathcal{U}_{r+n-1}^+ \otimes_{\mathbb{F}_p}
{\rm Fr}'^{r+n-1}
\left( \mathcal{U}_1^+ \right)}   @>>>  
{\mathcal{U}^+}
  \end{CD}
\] 
Here the upper and the left maps are 
$\mathbb{F}_p$-linear isomorphisms. 
Since the lower map gives an $\mathbb{F}_p$-linear isomorphism onto 
$\mathcal{U}_{r+n}^+$, so does the right map. Therefore, the result for 
the map $\mathcal{U}_r^+ \otimes_{\mathbb{F}_p} {\rm Fr}'^r
\left( \mathcal{U}_n^+\right) \rightarrow \mathcal{U}_{r+n}^+$ follows.  

The result for the second map in the theorem  
follows from that for the first one with $n=1$ and from induction on $r$. $\square$\\

Theorem \ref{mainthm1} immediately implies the following result. \\

\begin{Cor}\label{maincor1} 
Let $r \in \mathbb{Z}_{>0}$. 
The multiplication on $\mathcal{U}^+$ induces the following 
two $\mathbb{F}_p$-linear isomorphisms:  
\[
\mathcal{U}_r^+ \otimes_{\mathbb{F}_p} 
{\rm Fr}'^r \left( \mathcal{U}^+ \right) \rightarrow \mathcal{U}^+,\ \ \ 
\bigotimes_{i \geq 0} {\rm Fr}'^i \left( \mathcal{U}_1^+ \right) \rightarrow \mathcal{U}^+.
\] 
\end{Cor}

\noindent {\bf Remark.} In Lusztig's book  
\cite[Proposition 35.4.2 (b)]{lusztigbook}, 
a similar result is stated in the case of a quantum group. \\

We note that by symmetry, analogous results to those in this section hold 
for $\mathcal{U}^{-}$ as well.

\section{Some linear isomorphisms for $\mathcal{U}$}

The aim in this section is to give several $\mathbb{F}_p$-linear isomorphisms induced by 
multiplication in $\mathcal{U}$.  As in Section 4, throughout this section, 
$r$ denotes a fixed positive integer. 

Before describing  the results, we shall observe some properties 
about the subalgebra $\mathcal{U}^0$. 

To begin with, we note that the $\mathbb{F}_p$-algebra $\mathcal{U}^0$ is commutative. 
So the multiplication in $\mathcal{U}^0$ 
induces some $\mathbb{F}_p$-algebra isomorphisms unlike 
the case of $\mathcal{U}^+$ or $\mathcal{U}^-$. \\

\begin{Prop}
Let $n \in \mathbb{Z}_{> 0}$. Then the multiplication on $\mathcal{U}^0$ induces 
the following four $\mathbb{F}_p$-algebra isomorphisms:  
\[
\mathcal{U}_r^0 \otimes_{\mathbb{F}_p} {\rm Fr}'^r\left( \mathcal{U}_n^0\right) 
\rightarrow \mathcal{U}_{r+n}^0,\ \ 
\bigotimes_{i=0}^{r-1} {\rm Fr}'^i \left( \mathcal{U}_1^0\right) 
\rightarrow \mathcal{U}_{r}^0,
\]
\[
\mathcal{U}_r^0 \otimes_{\mathbb{F}_p} 
{\rm Fr}'^r \left( \mathcal{U}^0 \right) \rightarrow \mathcal{U}^0,\ \  
\bigotimes_{i \geq 0} {\rm Fr}'^i \left( \mathcal{U}_1^0 \right) \rightarrow \mathcal{U}^0. 
\]
\end{Prop}

\noindent {\itshape Proof.} It is enough to prove the results only for the first map 
$\mathcal{U}_r^0 \otimes_{\mathbb{F}_p} {\rm Fr}'^r\left( \mathcal{U}_n^0\right) 
\rightarrow \mathcal{U}_{r+n}^0$ and the second map 
$\bigotimes_{i=0}^{r-1} {\rm Fr}'^i \left( \mathcal{U}_1^0\right) 
\rightarrow \mathcal{U}_{r}^0$. 
It is easy to check that these maps are  $\mathbb{F}_p$-algebra homomorphisms, 
using the facts that the algebra $\mathcal{U}^0$ is commutative 
and that the map ${\rm Fr}' : \mathcal{U}^0 \rightarrow \mathcal{U}^0$ is an 
$\mathbb{F}_p$-algebra endomorphism. Moreover, 
by the remark just after Proposition \ref{binomial1}, 
we see that the first map takes the elements 
$
\left( \prod_{i=1}^l {h_i \choose m_i} \right) \otimes 
{\rm Fr}'^r \left( \prod_{i=1}^l {h_i \choose t_i} \right)
$
with $m_i \in \mathcal{N}_r$ and $t_i \in \mathcal{N}_n$ which form an 
$\mathbb{F}_p$-basis 
of $\mathcal{U}_r^0 \otimes_{\mathbb{F}_p} {\rm Fr}'^r\left( \mathcal{U}_n^0\right)$ 
to the elements $\prod_{i=1}^{l} {h_i \choose m_i +p^r t_i}$ which form an 
$\mathbb{F}_p$-basis of $\mathcal{U}_{r+n}^0$. Therefore, the result 
for the first map follows. 
The result for the second map follows from that for the first one and induction on 
$r$. $\square$ \\ \\

Now we introduce primitive idempotents in $\mathcal{U}_n^0$ for 
$n \in \mathbb{Z}_{>0}$. For $1 \leq i \leq l$, $n \in \mathbb{Z}_{>0}$, 
and $j \in \mathbb{Z}$, set 
\[
\mu_{i,j}^{(n)} = {h_i -j-1 \choose p^n-1} = 
\sum_{k=0}^{p^n-1} {-j-1 \choose p^n-1-k} {h_i \choose k}  
\]
in $\mathcal{U}^0$. This element lies in $\mathcal{U}_{n}^0$. Then the following hold: \\ 

\noindent $\bullet$ We have 
${h_i \choose t} \mu_{i,j}^{(n)} = {j \choose t} \mu_{i,j}^{(n)}$ 
for each $t \in \mathcal{N}_n$. 
\\

\noindent $\bullet$ The elements $\mu_{i,j}^{(n)}$ with $j \in \mathcal{N}_n$ are 
pairwise orthogonal idempotents in $\mathcal{U}_{n}^0$ satisfying 
$\sum_{j=0}^{p^n-1} \mu_{i,j}^{(n)} =1$. \\

\noindent $\bullet$ For $j,j' \in \mathbb{Z}$, we have 
\[ \mu_{i,j}^{(n)} = \mu_{i,j'}^{(n)} \Longleftrightarrow 
j \equiv j'\ ({\rm mod}\ p^n).\]

\noindent $\bullet$ Let $m,n \in \mathbb{Z}_{>0}$. Then for $s \in \mathcal{N}_m$ and 
$t \in \mathbb{Z}$, we have 
$\mu_{i,s+p^m t}^{(n+m)}=\mu_{i,s}^{(m)} {\rm Fr}'^m
\left( \mu_{i,t}^{(n)} \right)$. \\

\noindent For details of these facts, see the results in \cite[\S 4]{gros-kaneda15} 
for the case of type ${\rm A}_1$. These results easily imply those of a general type. 
For $n \in \mathbb{Z}_{> 0}$ and $\lambda \in X(T)$, set 
\[
\mu_{\lambda}^{(n)} = \prod_{k=1}^{l} \mu_{k, \langle \lambda, \alpha_k^{\vee} \rangle}^{(n)}
= \prod_{k=1}^{l} {h_k -\langle \lambda, \alpha_k^{\vee} \rangle -1 \choose p^n-1} 
\]
in $\mathcal{U}^0$. This element lies in $\mathcal{U}_{n}^0$. 
For $m \in \mathbb{Z}_{> 0}$, set 
\[
X_{m}(T)=\{ \lambda \in X(T)\ |\ \langle \lambda, \alpha_i^{\vee} \rangle 
\in \mathcal{N}_m, 
\forall i \in \{1, \dots, l\} \}.
\] 
Since $G$ is simply connected and simple, $X_{m}(T)$ is a system of representatives for 
$X(T)/p^m X(T)$ (see \cite[II, 3.15]{jantzenbook}). 
Using the above facts, the following proposition can be easily proved. \\

\begin{Prop}\label{idempotents}
Let $n \in \mathbb{Z}_{>0}$. The following hold. \\

\noindent {\rm (i)} For  $\lambda \in X(T)$, 
we have ${h_i \choose t} \mu_{\lambda}^{(n)} = 
{\langle \lambda, \alpha_i^{\vee} \rangle \choose t}  \mu_{\lambda}^{(n)}$ 
for each $i$ with $1 \leq i \leq l$ and each $t \in \mathcal{N}_n$. \\

\noindent {\rm (ii)} The elements $\mu_{\lambda}^{(n)}$ with $\lambda \in X_n(T)$ are 
pairwise orthogonal primitive idempotents in $\mathcal{U}_{n}^0$  satisfying 
$\sum_{\lambda \in X_n(T)} \mu_{\lambda}^{(n)} =1$ and form an $\mathbb{F}_p$-basis of 
$\mathcal{U}_{n}^0$. \\

\noindent {\rm (iii)} For  $\lambda,\lambda' \in X(T)$, we have 
\[ \mu_{\lambda}^{(n)} = \mu_{\lambda'}^{(n)} \Longleftrightarrow 
\lambda \equiv \lambda'\ ({\rm mod}\ p^nX(T)).\]

\noindent {\rm (iv)} Let $m \in \mathbb{Z}_{>0}$, $\lambda \in X_m(T)$, and 
$\lambda' \in X(T)$. Then we have 
$\mu_{\lambda +p^m \lambda'}^{(n+m)} = \mu_{\lambda}^{(m)} 
{\rm Fr}'^m \left( \mu_{\lambda'}^{(n)} \right)$. \\

\noindent {\rm (v)} We have 
$e_{\alpha}^{(m)} \mu_{\lambda}^{(n)} = 
\mu_{\lambda+m \alpha}^{(n)} e_{\alpha}^{(m)} $ for  $\lambda \in X(T)$,  
$\alpha \in \Phi$, and $m \in \mathbb{Z}_{\geq 0}$. 
\end{Prop}   
\

Now we shall introduce the notion of homogeneous elements in 
$\mathcal{U}_{\mathbb{Z}}$ or $\mathcal{U}$. 
Let $\gamma$ be a fixed element in $X(T)$
which is a $\mathbb{Z}_{\geq 0}$-linear combination of simple roots. If a 
nonzero  element 
$z \in \mathcal{U}_{\mathbb{Z}}^+$ (resp. $z \in \mathcal{U}_{\mathbb{Z}}^-$) is 
a $\mathbb{Z}$-linear combination of elements of the form 
\[
\prod_{i=1}^{t} e_{\gamma_i}^{(a_i)}\ \ \ (\mbox{resp. }\prod_{i=1}^{t} e_{-\gamma_i}^{(a_i)})
\]
satisfying $\gamma_i \in \Phi^+$, $a_i \in \mathbb{Z}_{> 0}$, and 
$\gamma = \sum_{i=1}^{t} a_i \gamma_i$, we say that $z$ is a homogeneous element 
and then denote $\gamma$ by $|z|$ (for example, see \cite[1.2.1]{lusztigbook}). 
Note also that 
any homogeneous element  $z \in \mathcal{U}_{\mathbb{Z}}^+$ 
(resp. $z \in \mathcal{U}_{\mathbb{Z}}^-$)  can be written as 
a $\mathbb{Z}$-linear combination of elements of the form 
${\bm e}^{({\bm a})}$ (resp. ${\bm f}^{({\bm a})}$) with 
${\bm a}=(a_1, \dots, a_{\nu}) \in (\mathbb{Z}_{\geq 0})^{\nu}$ satisfying 
$|z| = \sum_{i=1}^{\nu} a_i \beta_i$. 
In turn, a homogeneous element $z \in \mathcal{U}$ is defined as a nonzero element 
which is the image of a homogeneous element $z_0 \in \mathcal{U}_{\mathbb{Z}}$ under 
the reduction modulo $p$. Then $|z|$ is defined as $|z_0|$. 
\\

\begin{Prop}\label{fact3}
Let $x \in \mathcal{U}_{\mathbb{Z}}^+$ and $y \in \mathcal{U}_{\mathbb{Z}}^-$ be 
homogeneous elements. Then $xy-yx$ can be written as a finite sum 
$\sum_{{\bm a},{\bm b}}{\bm f}^{(\bm a)} z_{{\bm a}, {\bm b}} {\bm e}^{(\bm b)}$ 
for some 
$ z_{{\bm a}, {\bm b}} \in \mathcal{U}_{\mathbb{Z}}^0$,  where  
${\bm a}, {\bm b} \in (\mathbb{Z}_{\geq 0})^{\nu}$ satisfy 
$|{\bm e}^{(\bm b)}|< |x|$ and $|x|-|y|=|{\bm e}^{(\bm b)}|-|{\bm f}^{(\bm a)}|$.
\end{Prop}

\noindent {\itshape Proof.} 
The finiteness of the sum is clear by the condition of ${\bm a}$ and ${\bm b}$. 
For simplicity, we denote $e_{\alpha_i}$ and $e_{-\alpha_i}$ by $e_i$ and $f_i$ 
respectively. 
Since the subring $\mathcal{U}_{\mathbb{Z}}^+$ (resp. $\mathcal{U}_{\mathbb{Z}}^-$) 
is generated by all $e_{i}^{(n)}$ (resp. $f_{i}^{(n)}$) with $i \in \{1, \dots,l\}$ and 
$n \in \mathbb{Z}_{\geq 0}$, we may assume that 
\[
x= \prod_{k=1}^{s} e_{i_k}^{(b_{i_k})} \mbox{ \ and \ } y= \prod_{k=1}^{t} f_{j_k}^{(a_{j_k})}
\]  
with $i_k,j_k \in \{1, \dots, l\}$ and $a_{j_k}, b_{i_k} \in \mathbb{Z}_{> 0}$ for some 
$s,t \in \mathbb{Z}_{\geq 0}$ without loss of generality. We proceed in two steps. 
First, we deal with the case 
$\sum_{k=1}^t a_{j_k}=1$ (i.e. $y=f_j$). \\

\noindent {\bf Step 1.} If $j \in \{1, \dots, l\}$, then $xf_j-f_j x$ can be 
written as a finite sum $\sum_{{\bm c}}z_{\bm c} {\bm e}^{(\bm c)}$ for some 
$z_{\bm c} \in \mathcal{U}_{\mathbb{Z}}^0$, where 
${\bm c} \in (\mathbb{Z}_{\geq 0})^{\nu}$
 satisfies  $|{\bm e}^{(\bm c)}|=|x|-\alpha_j$. \\

We proceed by induction on $s$. It is clear when $s=0$ 
(i.e. $x=1$). So suppose that $s \geq 1$. 
We write $x=x' e_{i_s}^{(b_{i_s})}$, where $x'=\prod_{k=1}^{s-1} e_{i_k}^{(b_{i_k})}$. By induction, 
we have $x'f_j-f_j x'=\sum_{{\bm b}'} z'_{{\bm b}'} {\bm e}^{({\bm b}')}$ for some  
$z'_{{\bm b}'} \in \mathcal{U}_{\mathbb{Z}}^0$, 
where ${\bm b}' \in (\mathbb{Z}_{\geq 0})^{\nu}$ satisfies  
$|{\bm e}^{({\bm b}')}|=|x'|-\alpha_j=|x|-b_{i_s} \alpha_{i_s} -\alpha_j$. 

Suppose that $i_s \neq j$. By Proposition \ref{basicformulas} (iv), we have 
\begin{align*}
xf_j -f_j x &= x' e_{i_s}^{(b_{i_s})} f_j -f_j x' e_{i_s}^{(b_{i_s})} = 
x'  f_j e_{i_s}^{(b_{i_s})} -f_j x' e_{i_s}^{(b_{i_s})} \\
&= (x'f_j-f_j x')e_{i_s}^{(b_{i_s})} = \sum_{{\bm b}'} z'_{{\bm b}'} {\bm e}^{({\bm b}')}e_{i_s}^{(b_{i_s})}
\end{align*}
and $|{\bm e}^{({\bm b}')}e_{i_s}^{(b_{i_s})}|=|x|-\alpha_j$, as required. 

Suppose that $i_s = j$. By Proposition \ref{basicformulas} (iii), 
we have $x'(h_j-b_j+1) = \theta x'$ for some 
$\theta \in \mathcal{U}_{\mathbb{Z}}^0$. Then 
by Proposition \ref{basicformulas} (ii), we have 
\begin{align*}
xf_j -f_j x &= x' e_{j}^{(b_{j})} f_j -f_j x' e_{j}^{(b_{j})} \\
&= x'  f_j e_{j}^{(b_{j})} +x'(h_j-b_j+1)e_j^{(b_j-1)}-f_j x' e_{j}^{(b_{j})} \\
&= \sum_{{\bm b}'} z'_{{\bm b}'} {\bm e}^{({\bm b}')}e_{j}^{(b_{j})}+\theta x' e_j^{(b_j-1)}
\end{align*}
and $|{\bm e}^{({\bm b}')}e_{j}^{(b_{j})}|=|x' e_j^{(b_j-1)}|=|x|-\alpha_j$, and hence Step 1 follows. 
\\

\noindent {\bf Step 2.} $xy-yx$ can be written as a finite sum 
$\sum_{{\bm a}, {\bm b}}{\bm f}^{(\bm a)} z_{{\bm a}, {\bm b}} {\bm e}^{(\bm b)}$ 
for some $ z_{{\bm a}, {\bm b}} \in \mathcal{U}_{\mathbb{Z}}^0$, where 
 ${\bm a}, {\bm b} \in (\mathbb{Z}_{\geq 0})^{\nu}$ satisfy  
$|{\bm e}^{(\bm b)}|< |x|$ and $|x|-|y|=|{\bm e}^{(\bm b)}|-|{\bm f}^{(\bm a)}|$. \\

We use induction on $\sum_{k=1}^{t} a_{j_k}$. 
Now that we are doing calculations only in $\mathcal{U}_{\mathbb{Z}}$, 
all the coefficients there are guaranteed to 
be integers. So we do not need to mind  
denominators appearing in the coefficients. We write 
$y=(1/a_{j_1}) f_{j_1} y'$, where $y'= f_{j_1}^{(a_{j_1}-1)} \prod_{k=2}^{t} f_{j_k}^{(a_{j_k})}$. 
Note that 
\[
xy-yx = \dfrac{1}{a_{j_1}} (xf_{j_1} y' - f_{j_1} y' x) 
= \dfrac{1}{a_{j_1}} (xf_{j_1}-f_{j_1} x)y' + \dfrac{1}{a_{j_1}} f_{j_1}(xy'-y'x).
\]
By induction, we can write 
$xy' -y'x = \sum_{{\bm c}, {\bm d}} {\bm f}^{({\bm c})} 
z'_{{\bm c}, {\bm d}} {\bm e}^{({\bm d})}$  
for some  $z'_{{\bm c}, {\bm d}} \in \mathcal{U}_{\mathbb{Z}}^0$, where 
${\bm c}, {\bm d} \in (\mathbb{Z}_{\geq 0})^{\nu}$ satisfy  
$|{\bm e}^{({\bm d})}| < |x|$ and  
$|x|-|y'|= |{\bm e}^{({\bm d})}|- |{\bm f}^{({\bm c})}|$.  Then we have 
\[
f_{j_1} (xy' -y'x) = \sum_{{\bm c}, {\bm d}} f_{j_1} {\bm f}^{({\bm c})} 
z'_{{\bm c}, {\bm d}} {\bm e}^{({\bm d})}
\]
and 
$|{\bm e}^{(\bm d)}|-|f_{j_1} {\bm f}^{(\bm c)}|= |x|-|y'|-\alpha_{j_1}=|x|-|y|$. 

On the other hand, consider the element $(xf_{j_1}-f_{j_1} x)y'$. By Step 1, we have 
$xf_{j_1}-f_{j_1} x= \sum_{\bm d} z''_{\bm d} {\bm e}^{({\bm d})}$ for some 
$z''_{\bm d} \in \mathcal{U}_{\mathbb{Z}}^0$, 
where ${\bm d} \in (\mathbb{Z}_{\geq 0})^{\nu}$ satisfies  
$|{\bm e}^{({\bm d})}| = |x|-\alpha_{j_1}$.  
So we fix such ${\bm d}$ with $z''_{\bm d} \neq 0 $ and consider the element 
$z''_{{\bm d}} {\bm e}^{({\bm d})} y'$. By induction, we have 
\[
 {\bm e}^{({\bm d})} y'- y' {\bm e}^{({\bm d})} =
 \sum_{{\bm c}', {\bm d}'}  {\bm f}^{({\bm c}')} z''_{{\bm c}', {\bm d}'} 
{\bm e}^{({\bm d}')}
\]
for some  $z''_{{\bm c}', {\bm d}'} \in \mathcal{U}_{\mathbb{Z}}^0$, where 
${\bm c}', {\bm d}' \in (\mathbb{Z}_{\geq 0})^{\nu}$ satisfy  
$|{\bm e}^{({\bm d}')}| < |{\bm e}^{({\bm d})}|$ and  
$|{\bm e}^{({\bm d})}|-|y'|= |{\bm e}^{({\bm d}')}|- |{\bm f}^{({\bm c}')}|$. Since both $y'$ and ${\bm f}^{({\bm c}')}$ are monomials, 
by Proposition \ref{basicformulas} (iii) there exist 
$\theta_{\bm d}, \theta_{{\bm c}', {\bm d}} \in \mathcal{U}_{\mathbb{Z}}^0$ such that 
$z''_{{\bm d}} y'=y' \theta_{\bm d}$ and 
$z''_{{\bm d}}{\bm f}^{({\bm c}')}={\bm f}^{({\bm c}')} \theta_{{\bm c}', {\bm d}}$. 
Then we have 
$z''_{{\bm d}} {\bm e}^{({\bm d})} y' =y' \theta_{\bm d} {\bm e}^{({\bm d})} +
\sum_{{\bm c}', {\bm d}'}  {\bm f}^{({\bm c}')} 
\theta_{{\bm c}', {\bm d}} z''_{{\bm c}', {\bm d}'} 
{\bm e}^{({\bm d}')}$. 
Moreover, we have $|{\bm e}^{({\bm d})}| < |x|$, 
$|{\bm e}^{({\bm d}')}| < |x|$, and 
$ |{\bm e}^{({\bm d})}|- |y'|=
 |{\bm e}^{({\bm d}')}|- |{\bm f}^{({\bm c}')}|=|x|-|y|$, 
and Step 2 is proved. $\square$ \\
\ 

By taking the reduction modulo $p$ in the previous proposition, we obtain the following. \\

\begin{Cor}\label{fact3'}
Let $x \in \mathcal{U}^+$ and $y \in \mathcal{U}^-$ be 
homogeneous elements. Then $xy-yx$ can be written as a finite sum 
$\sum_{{\bm a}, {\bm b}}{\bm f}^{(\bm a)} z_{{\bm a}, {\bm b}} {\bm e}^{(\bm b)}$ 
for some  $ z_{{\bm a}, {\bm b}} \in \mathcal{U}^0$, where
${\bm a}, {\bm b} \in (\mathbb{Z}_{\geq 0})^{\nu}$  satisfy 
$|{\bm e}^{(\bm b)}|< |x|$ and $|x|-|y|=|{\bm e}^{(\bm b)}|-|{\bm f}^{(\bm a)}|$.
\end{Cor}
\ 

Now we are ready to give main results in this section. \\

\begin{The}\label{mainthm2}
Let $n \in \mathbb{Z}_{> 0}$. Then  the multiplication on $\mathcal{U}$ induces  
the following two $\mathbb{F}_p$-linear isomorphisms:  
\[\mathcal{U}_r \otimes_{\mathbb{F}_p} {\rm Fr}'^r(\mathcal{U}_n) 
\rightarrow \mathcal{U}_{r+n},\ \ \ \bigotimes_{i=0}^{r-1} {\rm Fr}'^i \left( \mathcal{U}_1\right) 
\rightarrow \mathcal{U}_{r}.\] 
\end{The}

\noindent {\itshape Proof.} 
We first prove the result for the first map 
$\mathcal{U}_r \otimes_{\mathbb{F}_p} {\rm Fr}'^r(\mathcal{U}_n) 
\rightarrow \mathcal{U}_{r+n}$. Note that both $\mathcal{U}_r \otimes_{\mathbb{F}_p} {\rm Fr}'^r(\mathcal{U}_n) $ and 
$\mathcal{U}_{r+n}$ have the same dimension $p^{(r+n)(l+2\nu)}$ over $\mathbb{F}_p$. 
Since the elements  
${\bm f}^{({\bm a})} \mu_{\lambda}^{(r)} {\bm e}^{({\bm b})}$ with 
${\bm a}, {\bm b} \in (\mathcal{N}_r)^{\nu}$ and $\lambda \in X_r(T)$ 
(resp. 
${\bm f}^{({\bm a}')} \mu_{\lambda'}^{(n)} {\bm e}^{({\bm b}')}$ with 
${\bm a}', {\bm b}' \in (\mathcal{N}_n)^{\nu}$ and $\lambda' \in X_n(T)$) form 
an $\mathbb{F}_p$-basis of  $\mathcal{U}_r$ (resp. $\mathcal{U}_n$), 
it is enough to show that 
the elements ${\bm f}^{({\bm a})} \mu_{\lambda}^{(r)} {\bm e}^{({\bm b})}
{\rm Fr}'^r \left( {\bm f}^{({\bm a}')} \mu_{\lambda'}^{(n)} {\bm e}^{({\bm b}')} \right)$ are 
linearly independent over $\mathbb{F}_p$. Suppose that 
\[
\sum_{({\bm a}, {\bm a}', {\bm b}, {\bm b}', \lambda, \lambda')}
\eta({\bm a}, {\bm a}', {\bm b}, {\bm b}', \lambda, \lambda') 
{\bm f}^{({\bm a})} \mu_{\lambda}^{(r)} {\bm e}^{({\bm b})}
{\rm Fr}'^r \left( {\bm f}^{({\bm a}')} \mu_{\lambda'}^{(n)} {\bm e}^{({\bm b}')} \right) =0 
\eqno{(**)}
\]
in $\mathcal{U}$, where ${\bm a}=(a_1, \dots, a_{\nu}), {\bm b}
=(b_1, \dots, b_{\nu}) \in  (\mathcal{N}_r)^{\nu}$, 
${\bm a}'=(a'_1, \dots, a'_{\nu}) , {\bm b}'=(b'_1, \dots, b'_{\nu}) 
\in  (\mathcal{N}_n)^{\nu}$, $\lambda \in X_r(T)$,  
$\lambda' \in X_n(T)$, and 
$\eta({\bm a}, {\bm a}', {\bm b}, {\bm b}', \lambda, \lambda') \in \mathbb{F}_p$. 
We need to show that $\eta({\bm a}, {\bm a}', {\bm b}, {\bm b}', \lambda, \lambda') =0$ 
for all such $6$-tuples $({\bm a}, {\bm a}', {\bm b}, {\bm b}', \lambda, \lambda')$. 
So suppose that there exists a $6$-tuple 
$({\bm a}, {\bm a}', {\bm b}, {\bm b}', \lambda, \lambda')$ such that 
$\eta({\bm a}, {\bm a}', {\bm b}, {\bm b}', \lambda, \lambda') \neq 0$. 
Consider the element 
${\bm f}^{({\bm a})} \mu_{\lambda}^{(r)} {\bm e}^{({\bm b})}
{\rm Fr}'^r \left( {\bm f}^{({\bm a}')} \mu_{\lambda'}^{(n)} {\bm e}^{({\bm b}')} \right) $ 
for a $6$-tuple $({\bm a}, {\bm a}', {\bm b}, {\bm b}', \lambda, \lambda')$ such that 
$\eta({\bm a}, {\bm a}', {\bm b}, {\bm b}', \lambda, \lambda') \neq 0$. 

Now we choose a unique element 
$\widetilde{\lambda} \in X_r(T)$ such that 
$\widetilde{\lambda} \equiv \lambda - \sum_{i=1}^{\nu} b_i \beta_i\ ({\rm mod}\ p^rX(T))$. 
Then by Proposition \ref{idempotents} (iii)-(v), we have 
\begin{align*}
\lefteqn{{\bm f}^{({\bm a})} \mu_{\lambda}^{(r)} {\bm e}^{({\bm b})}
{\rm Fr}'^r \left( {\bm f}^{({\bm a}')} \mu_{\lambda'}^{(n)} {\bm e}^{({\bm b}')} \right) } \\
&= {\bm f}^{({\bm a})} \mu_{\lambda}^{(r)} {\bm e}^{({\bm b})}
 {\rm Fr}'^r \left( {\bm f}^{({\bm a}')} \right) 
 {\rm Fr}'^r \left(  \mu_{\lambda'}^{(n)}  \right) 
{\rm Fr}'^r \left(  {\bm e}^{({\bm b}')} \right) \\
&= {\bm f}^{({\bm a})}  
{\bm e}^{({\bm b})} {\rm Fr}'^r \left( {\bm f}^{({\bm a}')} \right) 
\mu_{\lambda-\sum_{i=1}^{\nu} b_i \beta_i}^{(r)} {\rm Fr}'^r \left(  \mu_{\lambda'}^{(n)}  \right) 
{\rm Fr}'^r \left(  {\bm e}^{({\bm b}')} \right) \\
&= {\bm f}^{({\bm a})}  
{\bm e}^{({\bm b})} {\rm Fr}'^r \left( {\bm f}^{({\bm a}')} \right) 
\mu_{\widetilde{\lambda}}^{(r)} {\rm Fr}'^r \left(  \mu_{\lambda'}^{(n)}  \right) 
{\rm Fr}'^r \left(  {\bm e}^{({\bm b}')} \right) \\
&= {\bm f}^{({\bm a})}
{\bm e}^{({\bm b})} {\rm Fr}'^r \left( {\bm f}^{({\bm a}')} \right) 
\mu_{\widetilde{\lambda}+p^r \lambda'}^{(n+r)}{\rm Fr}'^r \left(  {\bm e}^{({\bm b}')} \right) 
\end{align*}
in $\mathcal{U}$. By Corollary \ref{fact3'}, we have 
\[
{\bm e}^{({\bm b})} {\rm Fr}'^r \left( {\bm f}^{({\bm a}')} \right)=
{\rm Fr}'^r \left( {\bm f}^{({\bm a}')} \right) {\bm e}^{({\bm b})} 
+\sum_{{\bm c},{\bm d}}{\bm f}^{({\bm c})} z_{{\bm c},{\bm d}} {\bm e}^{({\bm d})} 
\]
for some $z_{{\bm c},{\bm d}} \in \mathcal{U}^0$, where 
${\bm c},{\bm d} \in (\mathbb{Z}_{\geq 0})^{\nu}$ satisfy 
$| {\bm e}^{({\bm d})} |
< | {\bm e}^{({\bm b})} |$ and 
$| {\bm e}^{({\bm d})} |-| {\bm f}^{({\bm c})} | =
| {\bm e}^{({\bm b})} |-p^r | {\bm f}^{({\bm a}')} |$. Then by Proposition \ref{idempotents} (v), 
 we have 
\begin{align*}
\lefteqn{{\bm f}^{({\bm a})} \mu_{\lambda}^{(r)} {\bm e}^{({\bm b})}
{\rm Fr}'^r \left( {\bm f}^{({\bm a}')} \mu_{\lambda'}^{(n)} {\bm e}^{({\bm b}')} \right) } \\
&= {\bm f}^{({\bm a})}
\left( {\rm Fr}'^r \left( {\bm f}^{({\bm a}')} \right) {\bm e}^{({\bm b})} 
+\sum_{{\bm c},{\bm d}}{\bm f}^{({\bm c})} z_{{\bm c},{\bm d}} {\bm e}^{({\bm d})} \right)
\mu_{\widetilde{\lambda}+ p^r \lambda'}^{(n+r)}
{\rm Fr}'^r \left(  {\bm e}^{({\bm b}')} \right)  \\
&=  {\bm f}^{({\bm a})} {\rm Fr}'^r \left( {\bm f}^{({\bm a}')} \right) 
\mu_{\widetilde{\lambda}+\sum_{i=1}^{\nu} b_i \beta_i+ p^r \lambda'}^{(n+r)} 
{\bm e}^{({\bm b})} {\rm Fr}'^r \left(  {\bm e}^{({\bm b}')} \right) \\
& \ \ + \sum_{{\bm c},{\bm d}} {\bm f}^{({\bm a})} {\bm f}^{({\bm c})} 
z_{{\bm c},{\bm d}} \ 
\mu_{\widetilde{\lambda}+\sum_{i=1}^{\nu} d_i \beta_i+ p^r \lambda'}^{(n+r)} 
{\bm e}^{({\bm d})} {\rm Fr}'^r \left(  {\bm e}^{({\bm b}')} \right) \tag{$***$}.  
\end{align*}

Let $\mathcal{X}$ be a set of all 6-tuples 
$({\bm a}, {\bm a}', {\bm b}, {\bm b}', \lambda, \lambda')$ with 
${\bm a}, {\bm b}
\in  (\mathcal{N}_r)^{\nu}$, 
${\bm a}', {\bm b}'
\in  (\mathcal{N}_n)^{\nu}$, $\lambda \in X_r(T)$, and  
$\lambda' \in X_n(T)$ such that 
$\eta({\bm a}, {\bm a}', {\bm b}, {\bm b}', \lambda, \lambda')\neq 0$ and that 
there are no 
$6$-tuples $({\bm c}, {\bm c}', {\bm d}, {\bm d}', \mu, \mu')$ with 
${\bm c}, {\bm d} \in  (\mathcal{N}_r)^{\nu}$  
${\bm c}', {\bm d}'\in  (\mathcal{N}_n)^{\nu}$, $\mu \in X_r(T)$, $\mu' \in X_n(T)$ 
satisfying  
$\eta({\bm c}, {\bm c}', {\bm d}, {\bm d}', \mu, \mu') \neq 0$ and  
$| {\bm e}^{({\bm d})} |+p^r | {\bm e}^{({\bm d}')} |
> | {\bm e}^{({\bm b})} |+p^r | {\bm e}^{({\bm b}')} | $.  
Now the equalities $(**)$ and $(***)$ together with linearly independence over 
$\mathbb{F}_p$ of 
$\{ {\bm f}^{({\bm a})} z'_{{\bm a},{\bm b}} {\bm e}^{({\bm b})} \ |\ 
{\bm a}, {\bm b} \in (\mathbb{Z}_{\geq 0})^{\nu} \}$ with 
$z'_{{\bm a},{\bm b}} \in \mathcal{U}^0-\{0\}$ being fixed  imply that 
\[
\sum_{({\bm a}, {\bm a}', {\bm b}, {\bm b}', \lambda, \lambda') \in \mathcal{X}} 
\eta({\bm a}, {\bm a}', {\bm b}, {\bm b}', \lambda, \lambda')
{\bm F}_{{\bm a}, {\bm a}'} {\bm H}_{{\bm b}, \lambda, \lambda'} {\bm E}_{{\bm b}, {\bm b}'} =0
\]
in $\mathcal{U}$, where 
${\bm F}_{{\bm a}, {\bm a}'}={\bm f}^{({\bm a})} {\rm Fr}'^r \left( {\bm f}^{({\bm a}')} \right) $, 
${\bm H}_{{\bm b}, \lambda, \lambda'}=\mu_{\widetilde{\lambda}+\sum_{i=1}^{\nu} 
b_i \beta_i+ p^r \lambda'}^{(n+r)}$, and 
${\bm E}_{{\bm b}, {\bm b}'} ={\bm e}^{({\bm b})} {\rm Fr}'^r 
\left(  {\bm e}^{({\bm b}')} \right)$.  Since the elements 
${\bm F}_{{\bm a}, {\bm a}'} {\bm H}_{{\bm b}, \lambda, \lambda'} {\bm E}_{{\bm b}, {\bm b}'} $ for 
different 6-tuples 
$({\bm a}, {\bm a}', {\bm b}, {\bm b}', \lambda, \lambda')$ in $\mathcal{X}$ are 
linearly independent over $\mathbb{F}_p$, all 
$\eta({\bm a}, {\bm a}', {\bm b}, {\bm b}', \lambda, \lambda')$ for 
$({\bm a}, {\bm a}', {\bm b}, {\bm b}', \lambda, \lambda') \in \mathcal{X}$ must be zero. 
That is contradiction. Therefore, the result for the first map in the theorem follows. 

The result for the second map in the theorem  
follows from that for the first one with $n=1$ and from induction on $r$. 
$\square$ \\
\ 

As in the previous section, Theorem \ref{mainthm2} implies the following. \\

\begin{Cor}\label{maincor2}
The multiplication on $\mathcal{U}$ induces the following 
two $\mathbb{F}_p$-linear isomorphisms:  
\[
\mathcal{U}_r \otimes_{\mathbb{F}_p} 
{\rm Fr}'^r \left( \mathcal{U} \right) \rightarrow \mathcal{U},\ \ \ 
\bigotimes_{i \geq 0} {\rm Fr}'^i \left( \mathcal{U}_1 \right) \rightarrow \mathcal{U}.
\] 
\end{Cor}
\ 

\noindent {\bf Remark.} By Theorem \ref{mainthm2}, we also see  
that the multiplication in $\mathcal{U}^{\geq 0}$ induces the following four 
$\mathbb{F}_p$-linear isomorphisms:
\[
\mathcal{U}_r^{\geq 0} \otimes_{\mathbb{F}_p} {\rm Fr}'^r
\left( \mathcal{U}_n^{\geq 0}\right) 
\rightarrow \mathcal{U}_{r+n}^{\geq 0},\ \ \ 
\bigotimes_{i=0}^{r-1} {\rm Fr}'^i \left( \mathcal{U}_1^{\geq 0}\right) 
\rightarrow \mathcal{U}_{r}^{\geq 0},
\]
\[
\mathcal{U}_r^{\geq 0} \otimes_{\mathbb{F}_p} 
{\rm Fr}'^r \left( \mathcal{U}^{\geq 0} \right) \rightarrow \mathcal{U}^{\geq 0},\ \ \  
\bigotimes_{i \geq 0} {\rm Fr}'^i \left( \mathcal{U}_1^{\geq 0} \right) \rightarrow 
\mathcal{U}^{\geq 0}. 
\]
Of course, similar results for $\mathcal{U}^{\leq 0}$ also hold. \\ \\
\

\noindent {\large {\bf Acknowledgments}} \\

The author would like to thank the referee for carefully
reading the manuscript and giving some helpful comments. 

This work was supported by JSPS KAKENHI Grant Number JP18K03203.

\end{document}